\pdfoutput=1
\documentclass[12pt]{amsart}
\usepackage{amssymb}
\usepackage{booktabs}
\usepackage{url}
\usepackage{hyphenat}
\usepackage{mathtools}
\usepackage{cancel} % Allows cancelling terms
\usepackage[utf8]{inputenc}
\usepackage{ifpdf}
\ifpdf
  \usepackage[pdftex,final]{graphicx}
  \usepackage[pdftex,lmargin=1in,rmargin=1in,tmargin=1in,bmargin=1in]{geometry}
  \usepackage[bookmarks=true, bookmarksopen=true,%
    bookmarksdepth=3,bookmarksopenlevel=2,%
    colorlinks=true,%
    linkcolor=blue,%
    citecolor=blue,%
    filecolor=blue,%
    menucolor=blue,%
    pagecolor=blue,%
    urlcolor=blue]{hyperref}
  \let\textalt\texorpdfstring
\else
  \usepackage[dvips,final]{graphicx}
  \usepackage[dvips,lmargin=1in,rmargin=1in,tmargin=1in,bmargin=1in]{geometry}
  \newcommand{\textalt}[2]{#1}
  \makeatletter
  \let\origref\ref
  \DeclareRobustCommand\ref{\@ifstar\origref\origref}
  \makeatother
\fi

\newread\testin

\graphicspath{{draws/}{mpdraws/}}
\makeatletter
\def\input@path{{}{draws/}}
\makeatother

\newcommand{\RR}{\mathbb R}

\newcommand{\ZZ}{\mathbb Z}

\newcommand{\FF}{\mathbb F}

\newcommand{\Image}{\mathrm{Im}}

\newcommand{\co}{\colon}

\newcommand{\bdy}{\partial}

\newcommand{\lbracket}{[}
\newcommand{\rbracket}{]}

\DeclareMathOperator{\spin}{spin}
\newcommand{\SpinC}{\spin^c}

\DeclareMathOperator{\gr}{gr}

\theoremstyle{plain}
\newtheorem{theorem}{Theorem}
\numberwithin{equation}{section}
\newtheorem{proposition}[equation]{Proposition}
\newtheorem{lemma}[equation]{Lemma}

\theoremstyle{definition}

\theoremstyle{remark}
\newtheorem{example}[equation]{Example}
\newtheorem{remark}[equation]{Remark}
\newcommand{\HFa}{\widehat {\mathit{HF}}}

\newcommand{\CFa}{\widehat {\mathit{CF}}}

\newcommand{\x}{\mathbf x}
\newcommand{\y}{\mathbf y}

\newcommand{\w}{\mathbf w}

\newcommand{\Ainf}{\mathcal A_\infty}
\newcommand{\Alg}{\mathcal{A}}

\newcommand{\Idem}{\mathcal{I}}
\newcommand{\alphas}{{\boldsymbol{\alpha}}}
\newcommand{\betas}{{\boldsymbol{\beta}}}

\newcommand{\CFD}{\mathit{CFD}}

\newcommand{\CFA}{\mathit{CFA}}

\newcommand{\CFDa}{\widehat{\CFD}}

\newcommand{\CFK}{\mathit{CFK}}

\newcommand{\CFKm}{\CFK^-}

\newcommand{\CFAa}{\widehat{\CFA}}

\newcommand\DTP{\mathop{\widetilde\otimes}\nolimits}

\newcommand\Gen{\mathfrak{S}}
\renewcommand{\S}{\Gen}

\newcommand{\Heegaard}{\mathcal{H}}
\newcommand{\HD}{\Heegaard}

\makeatletter
\newcommand\honestalg[3]{\bigl\lbracket
\begin{smallmatrix} #1\@ifempty{#3}{}{&#3} \\ #2 \end{smallmatrix}
\bigr\rbracket}
\makeatother
\newcommand{\XX}{\mathbb{X}}
\newcommand{\OO}{\mathbb{O}}
\newcommand{\Xperm}{\sigma_\XX}
\newcommand{\Operm}{\sigma_\OO}
\newcommand{\Talpha}{\overline{\alpha}}
\newcommand{\Tbeta}{\overline{\beta}}
\newcommand{\Talphas}{\overline{\boldsymbol{\alpha}}}
\newcommand{\Tbetas}{\overline{\boldsymbol{\beta}}}
\newcommand{\CFP}{\mathit{CP}}

\newcommand{\CFPm}{\CFP^-}
\newcommand{\CFPDm}{\mathit{CPD}^-}
\newcommand{\CFPAm}{\mathit{CPA}^-}
\newcommand{\CFPDAm}{\mathit{CPDA}^-}
\newcommand{\CFPDDm}{\mathit{CPDD}^-}
\newcommand{\CFPAAm}{\mathit{CPAA}^-}
\newcommand{\Torus}{T}
\newcommand{\cI}{\mathcal{I}}
\newcommand{\Category}{\mathcal{C}}
\newcommand{\ModCat}{\mathsf{Mod}}

\newcommand{\DerBounded}{\mathcal{D}^b}
\newcommand{\Ring}{\mathbb{A}}

\DeclareMathOperator{\Rect}{Rect}
\newcommand{\TRect}{\overline{\Rect}}
\newcommand{\TRectEmpt}{\TRect^\circ}
\newcommand{\RectEmpt}{\Rect^\circ}
\DeclareMathOperator{\Half}{Half}
\newcommand{\HalfEmpt}{\Half^\circ}
\DeclareMathOperator{\Strip}{Strip}
\newcommand{\StripEmpt}{\Strip^\circ}

\newcommand{\Braids}{\mathcal{B}}

\DeclareMathOperator{\Cross}{cr}
\DeclareMathOperator{\smooth}{smooth}

\newcommand\Perm[1]{\lbracket\begin{array}{@{}*{10}{c@{\,}}}#1\end{array}\!\rbracket}

\begin{document}
\title[Slicing planar grid diagrams]{Slicing planar grid diagrams: a
  gentle introduction to bordered Heegaard Floer homology}

\author[Lipshitz]{Robert Lipshitz}
\thanks{RL was supported by an NSF Mathematical Sciences Postdoctoral
  Research Fellowship.}
\address{Department of Mathematics, Columbia University\\
  New York, NY 10027}
\email{lipshitz@math.columbia.edu}

\author[Ozsv\'ath]{Peter Ozsv\'ath}
\thanks{PO was supported by NSF grants number DMS-0505811 and FRG-0244663.}
\address {Department of Mathematics, Columbia University\\ New York, NY 10027}
\email {petero@math.columbia.edu}

\author[Thurston]{Dylan P.~Thurston}
\thanks {DPT was supported by a Sloan Research Fellowship.}
\address{Department of Mathematics,
         Barnard College,
         Columbia University\\
         New York, NY 10027}
\email{dthurston@barnard.edu}

\date{August 26, 2008}

\begin{abstract}
  We describe some of the algebra underlying the decomposition of
  planar grid diagrams. This provides a useful toy model for an
  extension of Heegaard Floer homology to 3-manifolds with
  parametrized boundary. This paper is meant to serve as a gentle
  introduction to the subject, and does not itself have immediate topological
  applications.
\end{abstract}

%\primaryclass{}
%\secondaryclass{}
%\keywords{}

\maketitle

\tableofcontents

\section{Introduction}
\label{sec:introduction}
The Heegaard Floer homology groups of Ozsv\'ath and Szab\'o are
defined in terms of holomorphic curves in Heegaard diagrams.
In~\cite{LOT1}, Heegaard Floer homology is extended to three-manifolds
with (parameterized) boundary, by studying holomorphic curves in
pieces of Heegaard diagrams. The resulting invariant, \emph{bordered
  Heegaard Floer homology}, has the following form.  To an oriented
surface $F$ (together with an appropriate Morse function on $F$),
bordered Heegaard Floer associates a differential graded algebra
$\Alg(F)$.  To a three-manifold $Y$ together with a homeomorphism
$F\to\bdy Y$, bordered Heegaard Floer associates a right ($\Ainf)$
module $\CFAa(Y)$ over $\Alg(F)$ and a left (differential graded)
module $\CFDa(Y)$ over $\Alg(-F)$. (Here, $-F$ denotes $F$ with its
orientation reversed.)  These modules, which are well-defined up to
homotopy equivalence, relate to the closed Heegaard Floer homology
group $\HFa$ via the following \emph{pairing theorem}:
\begin{theorem}[\cite{LOT1}] Suppose that $Y=Y_1\cup_F Y_2$. Then $\CFa(Y)\simeq \CFAa(Y_1)\DTP_{\Alg(F)}\CFDa(Y_2)$.
\end{theorem} 
(Recall that $\CFa(Y)$ is the chain complex underlying the Floer homology
group $\HFa(Y)$. The notation $\DTP$ denotes the derived tensor product,
and the symbol $\simeq$ denotes quasi-isomorphism.)

The definitions of the invariants $\CFAa$ and $\CFDa$ are,
unfortunately, somewhat involved. There are two kinds of complications
which obscure the basic ideas involved:
\begin{itemize}
\item {\bf Analytic complications.} The definitions of the invariants
  $\CFAa$ and $\CFDa$ involve counting pseudo-holomorphic curves. In spite of
  much progress over the last decades, holomorphic curve techniques
  remain somewhat technical, and often require seemingly unnatural
  contortions. To make matters worse, the analytic set up is, by
  necessity, somewhat nonstandard; in particular, it involves counting
  curves in a manifold with ``two kinds of infinities.''
\item {\bf Algebraic complications.} The invariant $\CFAa$ is, in general,
  not an honest module but only an $\Ainf$-module. While the subject
  of $\Ainf$ algebra is increasingly mainstream, it still adds a layer
  of obfuscation to the study of bordered Heegaard Floer
  homology. Further exacerbating the situation is a somewhat novel
  kind of grading.
\end{itemize}

In developing bordered Heegaard Floer homology we found it useful to
study a toy model, in terms of planar grid diagrams, in which these
complications are absent. It is the aim of the present paper to
present this toy model. We hope that doing so will make the definition
of bordered Heegaard Floer homology in \cite{LOT1} more palatable.

We emphasize up front that the main objects of study in this paper do
\emph{not} give topological invariants. Still, the algebra involved is
reminiscent of well-known objects from representation theory---in
particular, the nilCoxeter algebra---so this paper may be of further interest.

Throughout this paper, $\FF$ will denote the field with two
elements and $\Ring$ will denote $\FF[U_1,\dots,U_N]$ (for whichever
$N$ is in play at the time).

\emph{Acknowledgements.} The first author thanks the organizers of the
G{\"o}kova Geometry\slash Topology Conference for inviting him to participate in
this stimulating event. He also thanks C. Douglas for interesting
conversations in the summer of 2006. The authors also thank
M. Khovanov for pointing out the relationship of the algebra $\Alg$
with the nilCoxeter algebra.

%%% Local Variables: 
%%% mode: latex
%%% TeX-master: "PlanarMain"
%%% End: 

\section[Background]{Background on knot Floer homology and grid diagrams}
\label{sec:background}
We start by recalling the combinatorial definition of
Manolescu-Ozsv\'ath-Sarkar \cite{MOS06:CombinatorialDescrip} of the
knot Floer homology groups. 
% (See also Manolescu-Ozsv\'ath-Szab\'o-Thurston,
% \cite{MOST07:CombinatorialLink}, for further elucidation of this
% description.)

\begin{figure}
%Font is 12 point.
\centering
\includegraphics{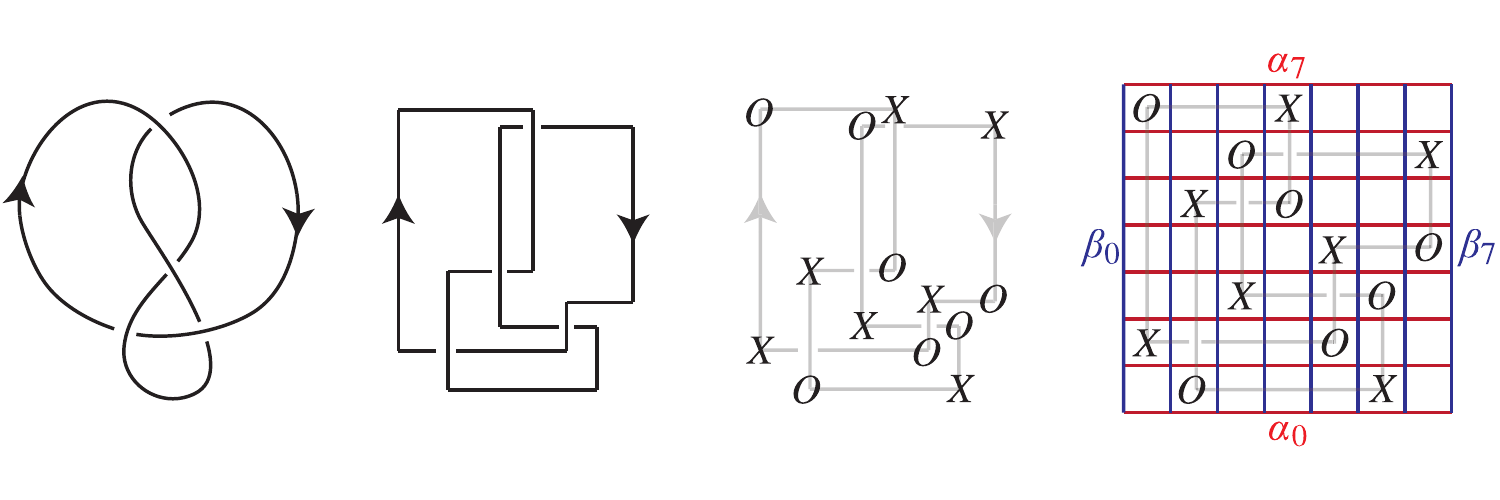}
\caption{\textbf{Representing a knot by a grid diagram.} Starting with
  a knot diagram $D$, one approximates $D$ using horizontal and
  vertical segments, so that crossings are always vertical over
  horizontal. Perturb the result so that no segments lie on the same
  horizontal or vertical line, and mark the endpoints alternately
  with $X$'s and $O$'s. The data of the knot is entirely encoded in
  these $X$'s and $O$'s, which we can see as sitting in the middle of
  squares on a piece of graph paper.}\label{fig:rep-by-grid}
\end{figure}

Let $K$ be an oriented knot in $S^3$. Choose a knot diagram $D$ for
$K$ such that
\begin{itemize}
\item $D$ is composed entirely of horizontal and vertical segments,
\item no two horizontal segments of $D$ have the same $y$-coordinate, and no
  two vertical segments of $D$ have the same $x$-coordinate, and
\item at each crossing, the vertical segment crosses over the
  horizontal segment.
\end{itemize}
(Every knot admits such a diagram; see Figure~\ref{fig:rep-by-grid}.)
The only data in such a diagram are the endpoints of the segments,
which we record by placing $X$'s and $O$'s at these endpoints,
alternately around the knot, and so that the knot is oriented from $X$
to $O$ along vertical segments. Notice that no two $X$'s (respectively
$O$'s) lie on the same horizontal or vertical line.

Let $\XX=\{X_i\}_{i=1}^N$ and $\OO=\{O_i\}_{i=1}^N$ denote the set of
$X$'s and $O$'s, respectively. Up to isotopy of the knot (and
renumbering of the $X_i$), we may assume that the coordinates of $X_i$
are $\left(i-\frac{1}{2},\Xperm(i)-\frac{1}{2}\right)$ for some
permutation $\Xperm\in S_N$. Then (after renumbering), the coordinates
of $O_i$ are $\left(i-\frac{1}{2},\Operm(i)-\frac{1}{2}\right)$ for some
permutation $\Operm\in S_N$. The data $(\RR^2,\XX,\OO)$ is a
\emph{planar grid diagram} for the knot $K$.

We can also view $\XX$ and $\OO$ as subsets of the torus
$\Torus=\RR^2/\langle (N,0),(0,N)\rangle.$ The data $(\Torus,\XX,\OO)$
is a \emph{toroidal grid diagram} for the knot $K$. It is easy to
recover the knot $K$ (up to isotopy) from the toroidal grid diagram
$(\Torus,\XX,\OO)$. We call the process of passing from a planar grid
diagram to a toroidal grid diagram \emph{wrapping}. The inverse
operation of passing from a toroidal grid diagram to a planar grid
diagram, which depends on a choice of two circles in $\Torus$, we call
\emph{unwrapping}.

The $N+1$ lines $\alpha_i=\{y=i\}\subset\RR^2$, $i=0,\dots,N$, descend
to $N$ disjoint circles $\Talpha_i$ in the torus $\Torus$, with
$\Talpha_0=\Talpha_N$. Similarly, the $N+1$ lines
$\beta_i=\{x=i\}\subset\RR^2$, $i=0,\dots,N$, descend to $N$ disjoint
circles $\Tbeta_i$ in $\Torus$. Notice that each $\alpha_j$
(respectively $\Talpha_j$) intersects each $\beta_i$ (respectively
$\Tbeta_j$) in a single point. Set $\alphas=\bigcup_{i=0}^N\alpha_i$,
$\Talphas=\bigcup_{i=1}^N\Talpha_i$, $\betas=\bigcup_{i=0}^N\beta_i$
and $\Tbetas=\bigcup_{i=1}^N\Tbeta_i$. We view the $\Talpha_i$ as
``horizontal'' and the $\Tbeta_i$ as ``vertical''. This means that
components of
$\Torus\setminus\left(\Talphas\cup\Tbetas\right)$ (little rectangles)
have, for instance,
\emph{lower left} corners, \emph{lower right} corners, and so on.

We define the knot Floer chain complex $\CFKm(K)$ as follows. Let
$\Ring=\FF[U_1,\dots,U_N]$. By a
\emph{toroidal generator} we mean an $N$-tuple of points
$\x=\{x_i\in\Talpha_{\sigma(i)}\cap\Tbeta_i\}$, one on each
$\Talpha$-circle and one on each $\Tbeta$-circle. Generators, then,
are in bijection with the permutation group $S_N$---but this bijection
depends on a choice of unwrapping.  Let $\S(\Torus,\XX,\OO)$ denote the set
of generators. The knot Floer complex $\CFKm(K)$ is freely generated
over $\Ring$ by $\S(\Torus,\XX,\OO)$.

For two generators $\x=\{x_i\}$ and $\y=\{y_i\}$, we define a set
$\TRect(\x,\y)$. The set $\TRect(\x,\y)$ is empty unless all but two
of the $x_i$ agree with corresponding $y_i$.  In that case, let $\{i,j\}
= \{k\mid x_k\neq y_k\}$; then
$\TRect(\x,\y)$ is the set of embedded rectangles $R$ in $\Torus$
with boundary on $\alphas\cup\betas$, and such that $x_i$ and $x_j$
are the lower-left and upper-right corners of~$R$ (in either order),
and $y_i$ and $y_j$ are the upper-left and lower-right corners of~$R$
(in either order). (Consequently,
$\TRect(\x,\y)$ always has either zero or two elements.) Call a
rectangle $R\in\TRect(\x,\y)$ \emph{empty} if the interior of $R$
contains no point in $\x$, and define $\TRectEmpt(\x,\y)$ to be the
set of empty rectangles in $\TRect(\x,\y)$. Given a rectangle $R$,
define $O_i(R)$ to be $1$ if $O_i$ lies in the interior of $R_i$ and
zero otherwise. Define $X_i(R)$ similarly, and set
$\OO(R)=\sum_{i=1}^NO_i(R)$ and $\XX(R)=\sum_{i=1}^NX_i(R)$.  Set
$U(R) = \prod_i U_i^{O_i(R)}$.

Now, define
\begin{equation}\label{eq:1}
  \bdy\x=\sum_{\y\in\S(\Torus,\XX,\OO)}
  \sum_{\begin{subarray}{c}R\in\TRectEmpt(\x,\y)\\
      \XX(R)=0 \end{subarray}}
  U(R)\cdot \y.
\end{equation}

\begin{lemma}
  Formula~(\ref{eq:1}) defines a differential, i.e., $\bdy^2=0$.
\end{lemma}
This is not hard to prove \cite[Proposition
2.8]{MOST07:CombinatorialLink}. See
Figure~\ref{fig:d-square-zero} for some of the cases.

By composing rectangles, we get more complicated regions in~$T$,
called \emph{domains}.  By a \emph{domain connecting $\x$ to
  $\y$} we mean a cellular two-chain $B$ in
$(\Torus,\Talphas\cup\Tbetas)$ with the following property. Let
$\bdy_{\alpha}B$ denote the intersection of $\bdy B$ with
$\alphas$. Then we require $\bdy(\bdy_{\alpha}B)=\y-\x$.  We can
define $O_i(B)$, $X_i(B)$, $\OO(B)$, $\XX(B)$, and $U(B)$ in the same
way as for rectangles.

There are two $\ZZ$-gradings on $\CFKm(K)$, the \emph{Maslov} or
\emph{homological} grading, denoted $\mu$, and the \emph{Alexander
  grading}, denoted $A$. These have the property that $\bdy$ preserves
$A$ and lowers $\mu$ by $1$. We give the combinatorial
characterization of $A$ and~$\mu$ from~\cite{MOST07:CombinatorialLink}, up
to an overall shift. First, some
notation. Given sets $E$ and $F$ in
$\RR^2$, let $\cI(E,F)$ denote the number of pairs $(e,f)\in E\times
F$ such that $e$ lies to the lower left of $f$ (i.e., the number of pairs
$e=(e_1,e_2)\in\RR^2$ and $f=(f_1,f_2)\in\RR^2$, such that $e_1<f_1$
and $e_2<f_2$).

Now, fix an unwrapping $(\RR^2,\XX,\OO)$ of the diagram
$(\Torus,\XX,\OO)$, so a generator $\x\in\S(\Torus,\XX,\OO)$
corresponds to a $N$-tuple of points $u(\x)$ in $\RR^2$. Then, for
some constants $C_A$ and $C_\mu$ depending on the diagram and the
unwrapping (but not on $\x$),
\begin{align*}
  A(\x)&=\cI(\XX,\x)-\cI(\OO,\x)+C_A\\
\mu(\x)&=\cI(\x,\x)-2\cI(\OO,\x)+C_\mu,
\end{align*}
cf.~\cite[Formulas (1) and (2)]{MOST07:CombinatorialLink}, bearing in
mind that $I(\XX,\x)$ differs from $I(\x,\XX)$ by a constant. Together
with the property that $A(U_i)=-1$ and $\mu(U_i)=-2$ this
characterizes $A$ and $\mu$ up to overall additive constants.

A fundamental result of
Manolescu-Ozsv\'ath-Sarkar~\cite{MOS06:CombinatorialDescrip} states
that the complex $\CFKm(K)$ defined above is bi-graded homotopy
equivalent to the complex $\CFKm(K)$ defined by
Ozsv\'ath and Szab\'o~\cite{OS05:HFL} and also
by Rasmussen~\cite{Rasmussen03:Knots}. It follows, in particular, that
the homotopy type of $\CFKm(K)$ is independent of the toroidal grid
diagram for $K$. The fact that the homotopy type of $\CFKm(K)$ depends
only on the knot $K$ can also be proved
combinatorially~\cite{MOST07:CombinatorialLink}.
\subsection{Planar Floer Homology}
In this paper we will study a modification of the grid diagram
construction of $\CFKm$, which we call the \emph{planar Floer
  homology} and denote $\CFPm$, obtained by replacing toroidal grid
diagrams by planar grid diagrams throughout the definition of
$\CFKm$. In the planar setting, when we have $N$ different $X$'s we will have
$N+1$ different $\alpha$- (respectively $\beta$-) lines: we view the process of
wrapping the diagram as identifying $\alpha_0$ with $\alpha_{N}$, and
$\beta_0$ with $\beta_{N}$. Thus, a generator over $\Ring$ of the
complex $\CFPm(\XX,\OO)$ is an $(N+1)$-tuple of points
$\x=\{x_i\in\alpha_{\sigma(i)}\cap\beta_i\}_{i=0}^N$. The set
$\S(\RR^2,\XX,\OO)$ is in canonical bijection with the symmetric group
$S_{N+1}$.

\begin{figure}
  \includegraphics{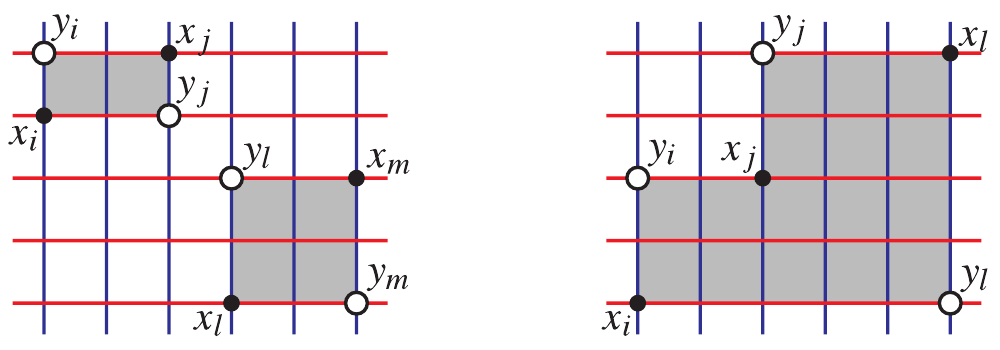}
  \caption{\textbf{Illustration of why $\bdy^2=0$ for planar Floer
      homology.} Left: The contributions to the coefficient of~$\y$
    from taking the two shaded rectangles in the two orders
    cancel. Right: This ``L''-shaped domain can be
    decomposed into two rectangles in two different ways, by making
    either a horizontal cut or a vertical cut. These two contributions
    cancel.}\label{fig:d-square-zero}
\end{figure}

Given generators $\x$ and $\y$ in $\S(\RR^2,\XX,\OO)$, let
$\RectEmpt(\x,\y)$ denote the set of empty rectangles in $\RR^2$
connecting $\x$ to $\y$; for each $\x$ and $\y$ the set
$\RectEmpt(\x,\y)$ is either empty or has a single element. The
differential on $\CFPm$ is defined analogously to
Formula~(\ref{eq:1}):
\begin{equation}\label{eq:2}
  \bdy\x=\sum_{\y\in\S(\RR^2,\XX,\OO)}
  \sum_{\begin{subarray}{c}R\in\RectEmpt(\x,\y)\\
      \XX(R)=0 \end{subarray}}
  U(R)\cdot \y.
\end{equation}
\begin{lemma}
  Formula~(\ref{eq:2}) defines a differential, i.e., $\bdy^2=0$.
\end{lemma}
The proof, which is a strict sub-proof of the proof for toroidal grid
diagrams, is illustrated in Figure~\ref{fig:d-square-zero}.

The complex $\CFPm(\XX,\OO)$ has Alexander and Maslov gradings $A$ and
$\mu$, defined exactly as they were for $\CFKm(K)$. We fix the
additive constants by setting
\begin{align*}
  A(\x)&=\cI(\XX,\x)-\cI(\OO,\x)\\
\mu(\x)&=\cI(\x,\x)-2\cI(\OO,\x).
\end{align*}

\textbf{\textit{Warning:} The homotopy type of the complex
  $\CFPm(\XX,\OO)$ is \textit{not} an invariant of the underlying
  knot $K$.} This is illustrated in
Example~\ref{example:not-invariant}. The results of this paper, thus,
do not directly give new topological invariants.

\begin{example}\label{example:not-invariant}
  Consider the planar grid diagrams for the unknot shown in
  Figure~\ref{fig:not-invariant}. The diagram on the left has
  $N=1$. The complex has two generators over $\FF[U_1]$, which we
  label with the permutations $\Perm{1&2}$ and $\Perm{2&1}$ in
  one-line notation.  (Here the one-line notation $\Perm{2&3&1}$, for
  instance, means the permutation $\{1\mapsto 2, 2\mapsto 3, 3\mapsto 1\}$.)
  The differential is trivial, so the
  homology of the complex is $\FF[U_1]^{\oplus 2}$.

  The diagram on the right has $N=2$. The complex has six
  generators. The differential is given by
  \begin{align*}
    \bdy\Perm{2&3&1}&=U_1\Perm{3&2&1}\\
    \bdy\Perm{3 & 1 & 2} &=U_2\Perm{3&2&1}\\
    \bdy\Perm{3&2&1}&=\bdy\Perm{1&2&3}=\bdy\Perm{1&3&2}=\bdy\Perm{2 &
      1 & 3}=0.
  \end{align*}
  The homology of the complex is 
  \[\FF\left\langle\Perm{3&2&1}\right\rangle
  \oplus \FF[U_1,U_2]\left\langle \Perm{1&2&3},\Perm{1&3&2},\Perm{2&1&3},U_2\Perm{2&3&1}+U_1\Perm{3&1&2}\right\rangle.
  \]
  This is certainly not the same as $\FF[U_1]^{\oplus 2}$.
  \begin{figure}
    %Font is 12 point.
    \centering
    \includegraphics{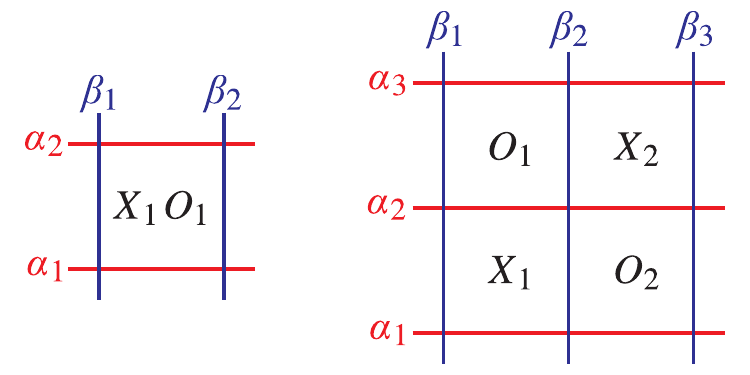}
    \caption{\textbf{Planar grid diagrams for the unknot.} Left: a
      diagram with $N=1$. Right: a diagram with $N=2$. The
      corresponding planar Floer complexes are not homotopy equivalent
      (Example~\ref{example:not-invariant}).}
    \label{fig:not-invariant}
  \end{figure}
\end{example}

% Local Variables: 
% mode: latex
% TeX-master: "PlanarMain"
% End: 

\section{Slicing planar grid diagrams}
\label{sec:slicing-planar-grid}
Fix a planar grid diagram $\HD=(\RR^2,\XX,\OO)$.  The goal of this
paper is to compute the complex $\CFPm(\HD)$ by cutting the diagram
vertically into pieces. (For now, we consider only cutting $\HD$ into
two pieces; we will consider more general cuttings in
Section~\ref{sec:freezing}.) We want to associate something
(ultimately, it will be a differential module) to each side, and
something (ultimately, it will be a differential graded algebra) to
the interface between the two sides. We want these to contain enough
information to reconstruct $\CFPm(\HD)$---but as little information as
possible beyond that, so as to be computable.

\begin{figure}
  \centering
  %Font is 12 point.
  \includegraphics{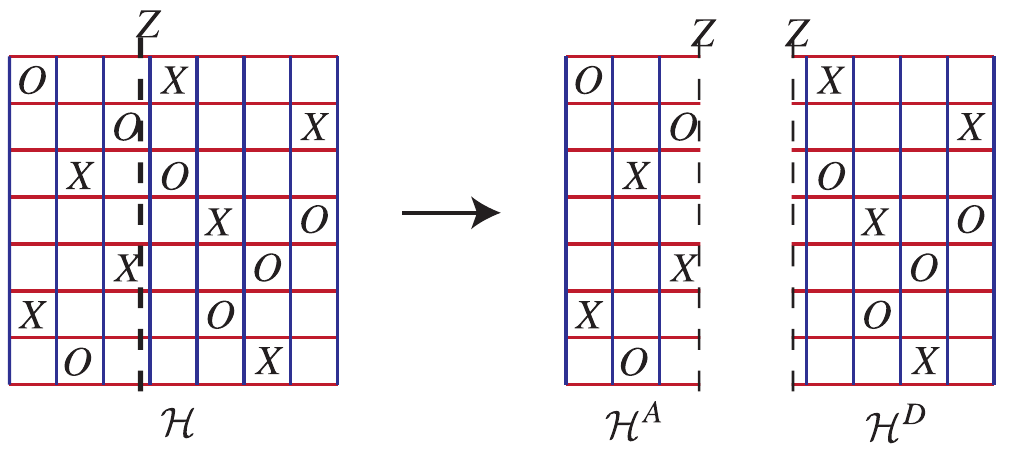}
  \caption{\textbf{Cutting a planar grid diagram.} The resulting
    diagrams $\HD^A$ and $\HD^D$ have $N=7$ and $k=3$.}
  \label{fig:cutting-diagram}
\end{figure}

So, let $Z$ be the vertical line $\{x=k-1/4\}$ and consider what each
side of $Z$ looks like. To the left of $Z$ we have $k$ vertical lines
$\beta_0,\dots,\beta_{k-1}$, as well as two injective maps
$\XX^A\co\{1,\dots,k\}\to\{1,\dots,N\}$ and
$\OO^A\co\{1,\dots,k\}\to\{1,\dots,N\}$. Similarly, to the right of
$Z$ we have $N+1-k$ vertical lines $\beta_k,\dots,\beta_{N}$, as well
as two injective maps $\XX^D\co\{k+1,\dots,N\}\to\{1,\dots,N\}$ and
$\OO^D\co\{k+1,\dots,N\}\to\{1,\dots,N\}$. There are also $N+1$
$\alpha$-lines, which intersect both sides of the diagram. Finally, at
the interface~$Z$ we see $N+1$ points $\{(i,k-1/4)\}_{i=0}^N$ where
the $\alpha_i$ intersect $Z$. See Figure~\ref{fig:cutting-diagram}.

Let $H^A$ denote the half-plane to the left of $Z$, and $H^D$ the
half-plane to the right of $Z$. We will call the data
$\HD^A=(H^A,\XX^A,\OO^A)$ or $\HD^D=(H^D,\XX^D,\OO^D)$ a \emph{partial
  planar grid diagram}. If we view $Z$ as oriented upwards then there
is a distinction between $\HD^A$ and $\HD^D$: for $\HD^A$ the induced
orientation of $Z$ agrees with the given one, while for $\HD^D$ the
induced orientation differs. We will call the first case ``type $A$''
and the second case ``type $D$.'' We say that $\HD^A$ has
\emph{height} $N+1$ and \emph{width} $k$, and $\HD^D$ has
\emph{height} $N+1$ and \emph{width} $N+1-k$.

Finally, a generator $\x=\{x_i\}_{i=0}^N$ corresponds to $k$ points
$\x^A=\{x_i\in\alpha_{\sigma^A(i)}\cap\beta_i\}_{i=0}^{k-1}$ to the left of $Z$
and $N+1-k$ points $\x^D=\{x_i\in\alpha_{\sigma^D(i)}\cap\beta_i\}_{i=k}^{N}$ to
the right of $Z$. Here, $\sigma^A$ is an injection
$\{0,\dots,k-1\}\to\{0,\dots,N\}$ and $\sigma^D$ is an injection
$\{k,\dots,N\}\to\{0,\dots,N\}$. For use later, let $\S(\HD^A)$ denote
the set of $k$-tuples
$\x^A=\{x_i\in\alpha_{\sigma^A(i)}\cap\beta_i\}_{i=0}^{k-1}$ where $\sigma^A$ is
an injection $\{0,\dots,k-1\}\to\{0,\dots,N\}$. Let $\S(\HD^D)$ denote
the set of $(N+1-k)$-tuples
$\x^D=\{x_i\in\alpha_{\sigma^D(i)}\cap\beta_i\}_{i=k}^{N}$ where $\sigma^D$ is an
injection $\{k,\dots,N\}\to\{0,\dots,N\}$.

% Local Variables: 
% mode: latex
% TeX-master: "PlanarMain"
% End: 

\section{Motivating the answer}
\label{sec:motivating-answer}
The purpose of this section is to motivate the answers which will
be described in later sections; thus, it can be skipped by the impatient
reader without sacrificing mathematical content.

We want to associate some kind of object, which with hindsight we
will call $\CFPAm(\HD^A)$ to $\HD^A$, and some other kind of object
$\CFPDm(\HD^D)$ to $\HD^D$. These should be objects in some
(algebraic) categories $\Category^A$ and $\Category^D$ associated to
the interface $Z$ (together perhaps with a little additional data). We
would also like a pairing map~$P$ from $\Category^A\times\Category^D$
to $\DerBounded(\Ring-\ModCat)$, the derived category of complexes over
the ground ring~$\Ring$, so that $\CFPm(\HD)=P(\HD^A,\HD^D)$. The
(derived) category
of chain complexes of (right/left) $\Alg$-modules for any
$\Ring$-algebra $\Alg$ admit such a pairing map, so this seems like a
reasonable example to keep in mind. (That is also how the story goes
in Khovanov homology \cite{Khovanov02:Tangles}, which is encouraging.)

Since a generator $\x$ of $\CFPm(\HD)$ decomposes as a pair
$(\x^A,\x^D)$, it seems reasonable that $\CFPAm(\HD^A)$ would be
generated---in some sense to be determined---by $\S(\HD^A)$ and that
$\CFPDm(\HD^D)$ would be generated by $\S(\HD^D)$.

Not every pair $(\x^A,\x^D)\in\S(\HD^A)\times\S(\HD^D)$ corresponds to
a generator in $\S(\HD)$: the necessary and sufficient condition is
that the images of the injections $\sigma^A$ and $\sigma^D$ be
disjoint. It seems reasonable that our putative $\Alg$ would remember
this---that if $\sigma^A_1$ and $\sigma^A_2$ have different images
then corresponding generators $\x_1^A$ and $\x_2^A$ would ``live
over'' different ``objects'' in $\Alg$. In the language of
differential graded categories (see, e.g.,
\cite{Keller06:DGCategories}), this makes sense; for algebras this
can be encoded via idempotents. That is, suppose $\Alg$ has
$\binom{N+1}{k}$ different primitive
idempotents $I_S$, one for each $k$-element subset $S$ of
$\{0,\dots,N\}$. Then we could say $\x^AI_S=\x^A$ if and only if
$S=\Image(\sigma^A)$, and $I_S\x^D=\x^D$ if and only if
$S\cap\Image(\sigma^D)=\emptyset$; otherwise these products
are~$0$. It then follows that an expression
of the form $\x^A\otimes_\Alg\x^D$ is nonzero if and only if
$(\x^A,\x^D)$ actually corresponds to a generator in $\S(\HD)$. We
will write $S(\x^A)$ to denote $\Image(\sigma^A)$, and $S(\x^D)$ to
denote $\{0,\dots,N\}\setminus\Image(\sigma^D)$.

There are three kinds of rectangles which contribute to the
differential on $\CFPm(\HD)$:
\begin{itemize}
\item Rectangles contained entirely in $\HD^A$. It seems reasonable
  that these should contribute to a differential on $\CFPAm(\HD^A)$,
  and there is an obvious way for them to do so.
\item Rectangles contained entirely in $\HD^D$. Again, it seems
  reasonable to let these contribute to a differential on
  $\CFPDm(\HD^D)$.
\item Rectangles which cross through the interface $Z$. It is somewhat
  less clear how to count these.
\end{itemize}

Let $R$ be a rectangle crossing through $Z$. Each of $\CFPAm(\HD^A)$
and $\CFPDm(\HD^D)$ see $Z$ as a half strip, and these half strips
should somehow be involved in the definitions of $\CFPAm(\HD^A)$ and
$\CFPDm(\HD^D)$. The rectangle $R$ intersects $Z$ in a segment running
from some $\alpha_i$ to some $\alpha_j$ (with $i < j$ by convention).
If $R$ is in $\Rect(\x,\y)$, with $\x=(\x^A,\x^D)$ and
$\y=(\y^A,\y^D)$, then the objects
(idempotents) associated to $\x^A$ and $\y^A$ differ:
$S(\y^A)=\left(S(\x^A)\setminus i\right)\cup j$. The objects $S(\x^D)$
and $S(\y^D)$ differ in the same way. So, we could view $R\cap Z$ as
an ``arrow'' from $S(\x^A)$ to $S(\y^A)$ or, in the
algebra language, as an element $\rho$ of $\Alg$ for which
$I_{S(\x^A)}\cdot \rho \cdot I_{S(\y^A)}=\rho$.

Actually, since a single rectangle in $\HD$ can be in $\Rect(\x,\y)$
for many different $\x$ and $\y$, the chord $R\cap Z$ gives many
arrows. More specifically, for any set $S$ with $i\in S$ and $j\notin
S$, $R\cap Z$ gives an arrow $\rho_{S,i,j}$, with the property that
$I_S\cdot\rho_{S,i,j}\cdot I_T=\rho_{S,i,j}$, where
$T=(S\setminus i)\cup j$. We
can view these as coming from a single element
$\rho_{i,j}=\sum_S\rho_{S,i,j}$ by multiplying with an idempotent. In
some sense, $\rho_{i,j}$ ``is'' $R\cap Z$.

With this in mind, there are two ways we can think of the effect of
the rectangle $R$ on one of the sides:
\begin{itemize}
\item It could start at $Z$, as the element $\rho_{i,j}$, and then
  come in to act on the module, moving one of the dots in the
  generator $\x$ to get the new generator~$\y$ (if not blocked). This
  is the point of view we will
  take for $\CFPAm$.
\item It could originate inside the partial diagram, and then propagate
  out to the boundary (if not blocked), leaving a residue $\rho_{i,j}$
  in $\Alg$ when it reaches the boundary. This is the point of view we
  will take for $\CFPDm$.
\end{itemize}
The two perspectives fit naturally with the pairing theorem: each
rectangle crossing the boundary starts in $\HD^D$, propagates out to
the boundary, and then propagates through to $\HD^A$.

More precisely, define $\CFPAm(\HD^A)$ to be generated \emph{over the
  base ring $\Ring$} by $\S(\HD^A)$. We have already defined an action
of the idempotents of $\Alg$ on $\CFPAm$. Define a right action of
$\Alg$ on $\CFPAm$ by setting $\x^A\cdot\rho_{i,j}=U(H)\cdot\y^A$ if
there is an empty half strip~$H$ connecting
$\x^A$ and $\y^A$ with rightmost edge equal to $\rho_{i,j}$ (and
not crossing any $X_k$). (Here $U(H)$ is the obvious extension of the
earlier notation to domains with boundary on~$Z$.) Define the product
to be zero
otherwise. Define the differential on $\CFPAm$ to count rectangles
entirely contained in $\HD^A$, in the obvious way.

Define $\CFPDm(\HD^D)$ to be ``freely'' generated as a left
$\Alg$-module by $\S(\HD^D)$. (More precisely, $\CFPDm$ is as free as
possible given the action of the idempotents we have already
defined. It is a direct sum of elementary modules, one for each
element of $\S(\HD^D)$.) Thus, the module structure on $\CFPDm$ is
rather dull. Define the differential on $\CFPDm$ as follows: given
generators $\x^D,\y^D\in\S(\HD^D)$, define
$\HalfEmpt(\rho_{i,j};\x^D,\y^D)$ to be the set of empty half strips
connecting $\x^D$ to $\y^D$ with boundary $\rho_{i,j}$; see
Figure~\ref{fig:DHalfStrip}. (The set
$\HalfEmpt(\rho_{i,j};\x^D,\y^D)$ is either empty or has a single
element.) Define
\[
  \bdy\x^D=\sum_{\y^D}\sum_{\begin{subarray}{c}R\in\RectEmpt(\x^D,\y^D)\\
    \XX(R)=0\end{subarray}} U(R)\cdot\y
  +\sum_{\y^D}\sum_{\rho_{i,j}}\sum_{\begin{subarray}{c}H\in\HalfEmpt(\rho_{i,j};\x^D,\y^D)\\
   \XX(H)=0\end{subarray}} U(H)\cdot\rho_{i,j}\y.
\]

\begin{remark} The $A$ in $\CFPAm$ is a mnemonic for the fact that the
  half-strips are included in the \emph{a}lgebra \emph{a}ction on
  $\CFPAm$. The $D$ in $\CFPDm$ is a mnemonic for the fact that the
  half-strips are included in the \emph{d}ifferential on $\CFPDm$.
\end{remark}

It is fairly clear that
$\CFPAm(\HD^A)\otimes_\Alg\CFPDm(\HD^D)=\CFPm(\HD)$. All rectangles
not crossing the interface are obviously accounted for.  If
$R\in\Rect(\x,\y)$ is a rectangle crossing the interface, with $R\cap
Z=\rho_{i,j}$, then
\[
\bdy(\x^A\otimes\x^D)=\x^A\otimes(\bdy\x^D)+\dots=\x^A\otimes\rho_{i,j}\cdot\y^D+\dots=\x^A\rho_{i,j}\otimes\y^D+\dots=\y^A\otimes\y^D+\cdots,
\]
as desired.

What is not clear---and, \emph{a priori}, not true---is that $\CFPAm$
and $\CFPDm$ are, in fact, chain complexes (differential modules) over
$\Alg$. Indeed, trying to make $\CFPDm$ into a module forces certain
relations---and a differential---on the algebra $\Alg$.

\begin{figure}
  \centering
  %Font is 12 point.
  \includegraphics{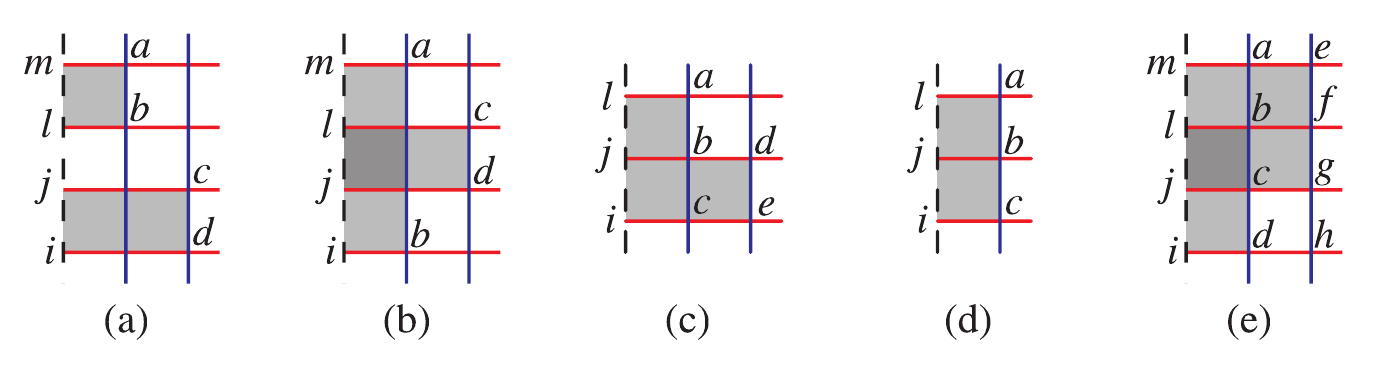}
  \caption{\textbf{Domains in $\HD^D$ forcing relations and a
      differential $\Alg$.} Part (a) forces $\rho_{i,j}$ and
    $\rho_{l,m}$ to commute. Part (b) forces $\rho_{i,m}$ and
    $\rho_{j,l}$ to commute. Part (c) forces
    $\rho_{i,j}\cdot\rho_{j,l}=\rho_{i,l}$. Part (d) forces the algebra to
    have a differential, and part (e) forces the product
    $\rho_{i,l}\cdot\rho_{j,m}$ to vanish.}
  \label{fig:D-forces}
\end{figure}

Consider the module $\CFPDm(\HD^D)$. In Part (a) of
Figure~\ref{fig:D-forces} is a plausible piece of $\HD^D$. One sees
here several generators; we single out $\{a,c\}$, $\{a,d\}$,
$\{b,c\}$ and $\{b,d\}$. Parts of the shaded region contribute to the
differential as follows:
\begin{align*}
\bdy\{a,c\}&=\rho_{l,m}\{b,c\}+\rho_{i,j}\{a,d\}+\cdots\\
\bdy\{b,c\}&=\rho_{i,j}\{b,d\}+\cdots\\
\bdy\{a,d\}&=\rho_{l,m}\{b,d\}+\cdots.
\end{align*}
(Here, the dots indicate contributions from regions of the diagram
other than the shaded one. The philosophy is that cancellation should
be local in $\HD^D$.)

Thus, one has
\[
\bdy^2\{a,c\}=\left(\rho_{l,m}\cdot\rho_{i,j}+\rho_{i,j}\cdot\rho_{l,m}\right)\{b,d\}+\cdots.
\]
So, in order to have $\bdy^2=0$ we should require that $\rho_{i,j}$
and $\rho_{l,m}$ commute.

Similarly, one sees by examining the shaded region in
Part (b) of Figure~\ref{fig:D-forces} that $\rho_{i,m}$ and
$\rho_{j,l}$ should commute.

In Part (c), consider the differentials
\begin{align*}
\bdy\{a,d\}&=\rho_{i,j}\{a,e\}+\rho_{i,l}\{c,d\}+\cdots\\
\bdy\{a,e\}&=\rho_{j,l}\{b,e\}+\cdots\\
\bdy\{c,d\}&=\{b,e\}+\cdots.
\end{align*}
Here,
\[
\bdy^2\{a,d\}=\rho_{i,j}\cdot\rho_{j,l}\{b,e\}+\rho_{i,l}\{b,e\}+\cdots.
\]
Thus, we should set $\rho_{i,j}\cdot\rho_{j,l}=\rho_{i,l}$---a relation
which looks rather reasonable in its own right.

Part (d) is a little trickier. Considering the generators
$\{a\}$, $\{b\}$ and $\{c\}$ we have
\begin{align*}
\bdy\{a\}&=\rho_{j,l}\{b\}+\rho_{i,l}\{c\}+\cdots\\
\bdy\{b\}&=\rho_{i,j}\{c\}+\cdots\\
\bdy\{c\}&=0+\cdots.
\end{align*}
Thus, it seems we have $\bdy^2\{a\}=\rho_{j,l}\cdot\rho_{i,j}\{c\}.$ One
might try setting $\rho_{j,l}\cdot\rho_{i,j}=0$, but it turns out this is
inconsistent with $\CFPAm$. Instead, we set (in this case)
\[
\bdy\rho_{i,l}=\rho_{j,l}\cdot\rho_{i,j}.
\]
Then it follows that $\bdy^2\{a\}=0$. Thus, we were forced to introduce
a differential on our algebra
$\Alg$.

Note that, in our example, $j\in S(\{a\})$. In general, we define
\[
\bdy(\rho_{S,i,l})=\sum_{\begin{subarray}{c}j\in S\\ i<j<l\end{subarray}}\rho_{S,j,l}\cdot\rho_{i,j}.
\]
This takes care of the example discussed above. The Leibniz rule
extends $\bdy$ to all of $\Alg$.

Part (e) is the most complicated. We will consider
$\bdy^2\{b,e\}$. We compute
\begin{align*}
\bdy\{b,e\}&=\{a,f\}+\rho_{j,l}\{c,e\}+\rho_{i,l}\{d,e\}\\
\bdy\{a,f\}&=\rho_{j,m}\{c,f\}+\rho_{i,m}\{d,f\}+\rho_{j,l}\{a,g\}+\rho_{i,l}\{a,h\}.\\
\bdy\left(\rho_{j,l}\{c,e\}\right)&=\rho_{j,l}\{a,g\}+\rho_{j,m}\{c,f\}+\rho_{j,l}\cdot\rho_{i,j}\{d,e\}\\
\bdy\left(\rho_{i,l}\{d,e\}\right)&=\rho_{j,l}\cdot\rho_{i,j}\{d,e\}+\rho_{i,l}\{a,h\}+\rho_{i,m}\{d,f\}+\rho_{i,l}\cdot\rho_{j,m}\{d,g\}.
\end{align*}
Most of the terms in $\bdy^2\{b,e\}$ cancel, but the term
$\rho_{i,l}\cdot\rho_{j,m}\{d,g\}$ does not. The offending domain is
shaded.

To resolve this difficulty, we impose the relation
$\rho_{i,l}\cdot\rho_{j,m}=0$ whenever $i<j<l<m$.

These are essentially all of the cases to check for $\CFPDm$; we will
verify this more carefully in Section~\ref{sec:type-D}.

Finally, consider the module $\CFPAm(\HD^A)$. One must check that the
relations we imposed on $\Alg$ are compatible with the action of
$\Alg$ on $\CFPAm(\HD^A)$; roughly, this follows by rotating the
pictures from
Figure~\ref{fig:D-forces} by $180$ degrees. We will discuss this more
thoroughly in Section~\ref{sec:type-A}.

These are the only relations we will need to impose on the algebra
$\Alg$. It turns out---we will see this next---that this algebra has a
clean description in terms of \emph{strand diagrams}.

%%% Local Variables: 
%%% mode: latex
%%% TeX-master: "PlanarMain"
%%% End: 

\section{The algebra associated to a slicing}
Fix integers $N+1$ and $k$, representing
the
height and width respectively of a partial planar grid diagram $\HD^A$. We will
define an algebra $\Alg_{N,k}$. We indicated, in a somewhat
roundabout manner, generators and relations for $\Alg_{N,k}$ in
Section~\ref{sec:motivating-answer}. We start by giving that
definition in a more orderly manner and then move on to a description
in terms of strand diagrams.

The algebra $\Alg_{N,k}$ is free as an $\Ring$-module. For each
$k$-element subset $S$ of $\{0,\dots,N\}$ there is a primitive
idempotent $I_S$, so that
\[
I_S\cdot I_T=
\begin{cases}
  I_S & \text{if $S=T$},\\
  0 & \text{otherwise.}
\end{cases}
\]

The algebra $\Alg_{N,k}$ is generated as an
$\Ring$-algebra by a set of elements $\rho_{S,i,j}$
(together with the idempotents). Here, $0\leq i<j\leq N$ and $S$ is a
$k$-element subset of $\{0,\dots,N\}$ such that $i\in S$ and $j\notin
S$. The relations with the idempotents are as follows:
\begin{align*}
I_T\cdot\rho_{S,i,j}&=
\begin{cases} \rho_{S,i,j} & \text{if $S=T$}\\
0 & \text{otherwise}
\end{cases}\\
\rho_{S,i,j}\cdot I_T&=
\begin{cases} \rho_{S,i,j} & \text{if $T=(S\setminus i)\cup j$}\\
0 & \text{otherwise}.
\end{cases}
\end{align*}
Set $\rho_{i,j}=\sum_S\rho_{S,i,j}$, so
$\rho_{S,i,j}=I_S\rho_{i,j}$. The relations we impose on $\Alg_{N,k}$ are:
\begin{alignat}{2}
  \rho_{i,j}\cdot\rho_{l,m}&=\rho_{l,m}\cdot\rho_{i,j}&\qquad & \text{for $j<l$ or $i<l<m<j$}\label{eq:rel1}\\
  \rho_{i,j}\cdot\rho_{l,m}&=0 &&\text{for $i<l<j<m$}\label{eq:rel2}\\
  \rho_{S,i,j}\cdot\rho_{j,l}&=\rho_{S,i,l}\label{eq:rel3} &&\text{for
    $j\notin S$.}
\end{alignat}
We also define a differential on $\Alg_{N,k}$ by setting
\[
\bdy(\rho_{S,i,j})=\sum_{\begin{subarray}{c}l\in S\\
    i<l<j\end{subarray}} \rho_{l,j}\cdot\rho_{i,l}
\]
and extending by the Leibniz rule.

Let $\Idem_{N,k}$ denote the subalgebra of $\Alg_{N,k}$ generated by
the idempotents.

We will check that $\bdy^2=0$ and that $\bdy$ has a consistent
extension to all of $\Alg_{N,k}$, but first we reinterpret this
algebra graphically, and introduce a grading.

\begin{figure}
  \centering
  \includegraphics{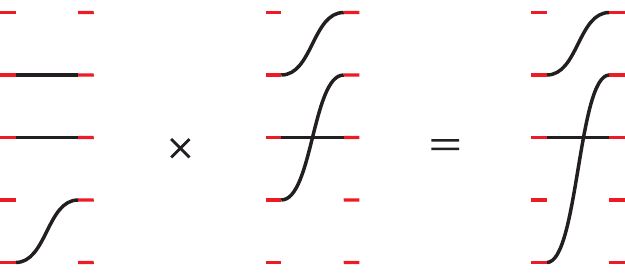}
  \caption{\textbf{The product on $\Alg_{4,3}$.} Two examples of
    upward-veering strand diagrams on $3$ strands and $5$ positions
    are shown left and center, and their product on the right.}
  \label{fig:up-flat-braid}
\end{figure}

Let $kI=\coprod_{i=1}^k[0,1]$, $\bdy_-kI=\coprod_{i=1}^k\{0\}$ and
$\bdy_+kI=\coprod_{i=1}^k\{1\}$. 
By an \emph{upward-veering strand diagram on $k$ strands and $N+1$
  positions} we mean a class $[\rho]$ of smooth maps 
\[
\rho\co (kI,\bdy_-kI,\bdy_+kI)\to
\left([0,1]\times[0,N],\{0\}\times\{0,\dots,N\},
  \{1\}\times\{0,\dots,N\}\right)
\]
such that $\rho'(t)\geq 0$ for all $t\in kI$, and such that the restrictions
$\rho|_{\bdy_-kI}$ and $\rho|_{\bdy_+kI}$ are injective, modulo
homotopy and reordering of the strands. (See
Figure~\ref{fig:up-flat-braid} for an illustration.) Let
$\Braids(N,k)$ denote the set of upward-veering strand diagrams on
$k$ strands and $N+1$ positions.

Given an element $[\rho]\in\Braids(N,k)$, let $\Cross([\rho])$ denote the
minimum number of crossings (double points) of any representative $\rho$ of $[\rho]$.

If $[\rho_1],[\rho_2]\in\Braids(N,k)$ are such that
$\bdy_+[\rho_1]=\bdy_-[\rho_2]$ then we can concatenate $\rho_1$ and
$\rho_2$ to obtain a new upward-veering strand diagram
$\rho_1\rho_2$. Note that $\Cross([\rho_1\rho_2])\leq\Cross([\rho_1])+\Cross([\rho_2])$. Let $\tilde{\Alg}_{N,k}$ denote the free
$\Ring$-module on $\Braids(N,k)$, and extend the
concatenation operation to a product on $\tilde{\Alg}_{N,k}$ by setting
\[
[\rho_1]\cdot[\rho_2]=
\begin{cases}
  [\rho_1\rho_2] & \text{if $\bdy_+[\rho_1]=\bdy_-[\rho_2]$ and
    $\Cross([\rho_1\rho_2])=\Cross([\rho_1])+\Cross([\rho_2])$}\\
  0 & \text{otherwise.}
\end{cases}
\]
This operation is obviously associative. The idempotents of
$\tilde{\Alg}_{N,k}$ are braids consisting of $k$ horizontal
strands, and as such are in bijection with the set of $k$-element
subsets of $\{1,\dots,N\}$.

\begin{figure}
  \centering
  %Font is 12 point.
  \includegraphics{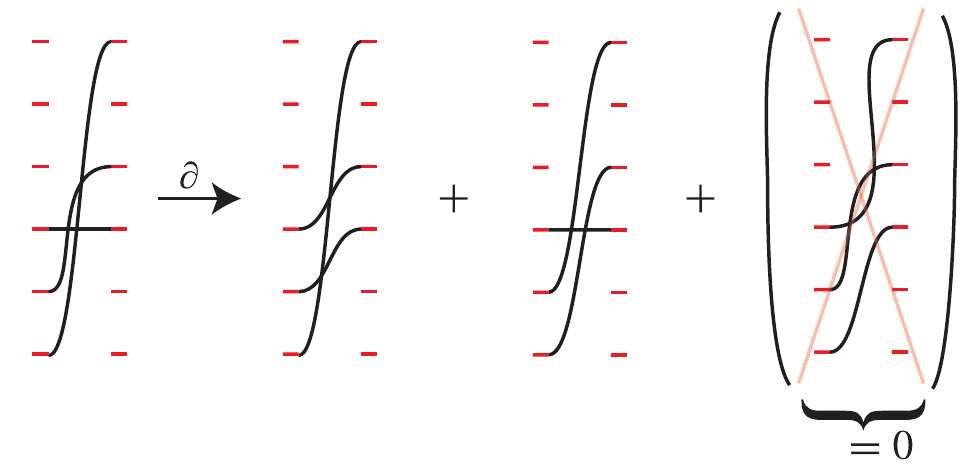}
  \caption{\textbf{The differential on $\Alg_{5,3}$.} Note that the
    term on the far right is not included in the differential of the
    term on the left because of the condition on the number of
    crossings $\Cross$.}
  \label{fig:ill-diff}
\end{figure}

We define a differential $\bdy$ on $\tilde{\Alg}_{N,k}$. Given
$[\rho]\in\Braids(N,k)$, with representative $\rho$, let
$\smooth(\rho)$ denote the multiset of strand diagrams obtained by
smoothing a single crossing in $\rho$. Then define
\[
\bdy[\rho]=\sum_{
  \begin{subarray}{c}
    \rho'\in\smooth(\rho)\\
    \Cross([\rho'])=\Cross([\rho])-1
  \end{subarray}
}[\rho'].
\]
See Figure~\ref{fig:ill-diff}.

\begin{figure}
  \centering
  %Font is 12 point
  \includegraphics{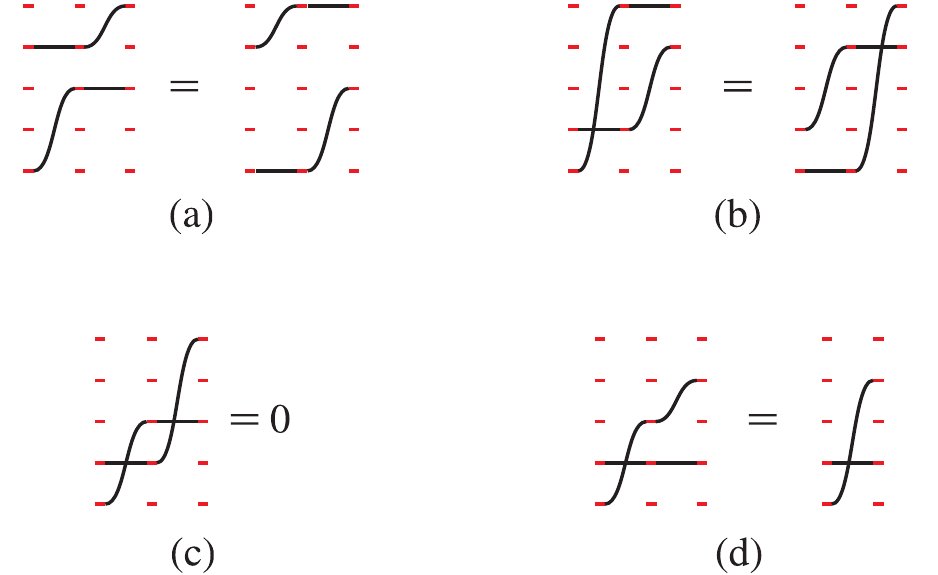}
  \caption{\textbf{The relations on $\Alg_{4,2}$.} Parts (a) and (b)
    correspond to relation~\eqref{eq:rel1}. Part (c) corresponds to
    relation~\eqref{eq:rel2}. Part (d) corresponds to
    relation~\eqref{eq:rel3}.}
  \label{fig:ill-relations}
\end{figure}

\begin{lemma}The algebra $\Alg_{N,k}$ is isomorphic to the algebra
  $\tilde{\Alg}_{N,k}$, via an isomorphism identifying the
  differentials.
\end{lemma}
\begin{proof}
  This is easy to check; see Figure~\ref{fig:ill-relations} for a
  convincing illustration that the relations agree. That the
  differentials agree is similarly straightforward.
\end{proof}

Provisionally, we define a grading on $\Alg_{N,k}$ by setting
$\gr([\rho])=\Cross([\rho])$.

\begin{proposition}The algebra $\Alg_{N,k}$ is a differential graded
  algebra. That is:
  \begin{enumerate}
  \item The differential satisfies $\bdy^2=0$.
  \item The differential satisfies the Leibniz rule $\bdy(ab)=(\bdy
    a)b+a(\bdy b)$.
  \item Multiplication has degree $0$.
  \item The differential has degree $-1$.
  \end{enumerate}
\end{proposition}
\begin{proof}
  All four parts are obvious from the description in terms of
  strand diagrams.
\end{proof}

\begin{remark} We have given two different definitions of
  $\Alg_{N,k}$. We could give a third, closely related to
  permutations: the algebra is generated over $\Ring$
  by bijective maps $f\co S\to T$ between $k$-element subsets of
  $\{0,\dots,N\}$, such that for all $i\in S$, $f(i)\geq i$. The
  function $\Cross$ is then the number of inversions of the map
  (i.e., the number of pairs of integers $i<j$ for which $f(i)>f(j)$), and
  the multiplication is composition if it is defined and preserves
  $\Cross$ and zero otherwise. See \cite[Section 3.1.1]{LOT1} for
  further discussion.
\end{remark}

The homological (Maslov) grading we want is not the same as $\gr$. In
fact, both the Maslov and Alexander gradings on $\Alg_{N,k}$ depend
not just on $N$ and $k$ but also on which rows contain $X$'s and $O$'s
to the left of $Z$.

More precisely, fix $k$-element subsets $L_X$ and $L_O$ of
$\{1/2,\dots,N-1/2\}$, which are  the $y$-coordinates of the
$X_i$'s and $O_i$'s contained in $\HD^A$ (including $X_k$). Given an
algebra element $a$, viewed as a strand diagram, let $L_X(a)$ denote
the intersection number of $a$ with the lines $y=\ell$ for $\ell\in
L_X$. (Equivalently, define $L_X(\rho_{i,j})=\#\{\ell\in L_X\mid i<\ell<j\}$
and extend to all of $\Alg_{N,k}$.) Define $L_O(a)$ similarly.

For $a\in\Alg_{N,k}$, define gradings $A$ and $\mu$ by
\begin{align*}
A(a)&=L_X(a)-L_O(a)\\
\mu(a)&=\Cross(a)-2L_O(a).
\end{align*}
It is clear that $A$ is preserved by multiplication and the
differential, and that multiplication preserves $\mu$ while the
differential drops $\mu$ by $1$.

% Local Variables: 
% mode: latex
% TeX-master: "PlanarMain"
% End: 

\section{The Type \textalt{$D$}{D} module}
\label{sec:type-D}
Fix a partial planar grid diagram $\HD^D$ of height $N+1$ and width
$N+1-k$. We will associate to $\HD^D$ a differential
$\Alg_{N,k}$-module.

We define a left action of the idempotents $\Idem_{N,k}$ on $\Ring\langle
\S(\HD^D)\rangle$, the free $\Ring$-module generated by the generators in
$\HD^D$ (see Section~\ref{sec:slicing-planar-grid}). Recall that a
generator $\x^D\in\S(\HD^D)$ corresponds to an injection
$\sigma_\x\co\{k,\dots,N+1\}\to\{0,\dots,N\}$. So, set
\[
I_S\x^D=
\begin{cases}
  \x^D&\text{if $S\cap\Image(\sigma_{\x^D})=\emptyset$}\\
  0 & \text{otherwise.}
\end{cases}
\]
As an $\Alg_{N,k}$-module, let
\[
\CFPDm(\HD^D)=\Alg_{N,k}\otimes_{\Idem_{N,k}}\Ring\langle
\S(\HD^D)\rangle.
\]
That is, the module $\CFPDm(\HD^D)$ is a direct sum of elementary
$\Alg_{N,k}$-modules, one for each generator in $\S(\HD^D)$.

\begin{figure}
  \centering
  \includegraphics{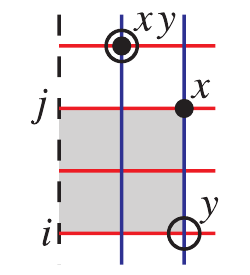}
  \caption{\textbf{An element (shaded) of $\Half(\rho_{i,j};\x,\y)$.}
    In fact, the region pictured lies in
    $\HalfEmpt(\rho_{i,j};\x,\y)$. It is also permitted for there to
    be some $O_i$ or $X_i$ in the domain (though in the latter case we
    will not, in fact, count the domain for the theory under
    discussion).}
  \label{fig:DHalfStrip}
\end{figure}

We next define the differential on $\CFPDm$. For generators $\x^D$
and $\y^D$, define $\RectEmpt(\x^D,\y^D)$ exactly as in
Section~\ref{sec:background}. Given generators $\x^D$, $\y^D$ and a
segment $\rho_{i,j}$ in $Z$, we define a set $\Half(\rho_{i,j};\x,\y)$, as
follows. Define $\Half(\rho_{i,j};\x,\y)$ to be empty unless
$x_i=y_i$ for all but one $i$. If $x_i=y_i$ for $i\neq j$ and the
$y$-coordinate of $x_j$ is (strictly) greater than the $y$-coordinate
of $y_j$, then let $\Half(\rho_{i,j};\x,\y)$ be the singleton set containing the
rectangle (or ``half-strip'') $H$ with upper right corner  $x_j$, and lower
right corner $y_j$, and left edge along the interface $Z$, where it is the segment
from $y=i$ to $y=j$. See Figure~\ref{fig:DHalfStrip}. Call a half
strip $H\in\Half(\rho_{i,j};\x^D,\y^D)$ \emph{empty} if the interior
of $H$ is disjoint from $\x^D$ (or equivalently from $\y^D$). Let
$\HalfEmpt(\rho_{i,j};\x^D,\y^D)$ denote the set of empty half strips
in $\Half(\rho_{i,j};\x^D,\y^D)$; this set has at most one element.

Now, for $\x^D$ a generator, define
\[
  \bdy\x^D=\sum_{\y^D\in\S(\HD^D)}
  \sum_{
    \begin{subarray}{c}
      R\in\RectEmpt(\x^D,\y^D)\\
      \XX(R)=0
    \end{subarray}
  }U(R)\cdot\y
  +\sum_{\y^D}\sum_{\rho_{i,j}}\sum_{
    \begin{subarray}{c}
      H\in\HalfEmpt(\rho_{i,j};\x^D,\y^D)\\
      \XX(R)=0
    \end{subarray}
  }U(H)\cdot\rho_{i,j}\y.
\]
We extend the definition via the Leibniz rule to all of $\CFPDm(\HD^D)$.

\begin{proposition}\label{prop:D-d-squared}The module $(\CFPDm,\bdy)$ is a differential
  module. That is, $\bdy^2=0$.
\end{proposition}
\begin{proof}
Since
\begin{align*}
\bdy^2(a\w^D)&=\bdy\left((\bdy a)\w^D+a(\bdy\w^D)\right)\\
&=\left(\bdy^2a\right)\w^D+2(\bdy a)(\bdy
\w^D)+a\left(\bdy^2\w^D\right)\\
&=a\left(\bdy^2\w^D\right),
\end{align*}
it suffices to show that the coefficient of $\y^D$ in $\bdy^2\w^D$ is
zero for any $\w^D,\y^D\in\S(\HD^D)$.

The remainder of the proof is similar to the combinatorial proof in
the closed case \cite[Proposition
2.8]{MOST07:CombinatorialLink}. Let $a_{\w^D,\x^D}$ denote the
coefficient of $\x^D$ in $\bdy\w^D$. Then the coefficient of $\y^D$ in
$\bdy^2\w^D$ is
\begin{equation}\label{eq:3}
\left(\sum_{\x^D} a_{\w^D,\x^D}\cdot a_{\x^D,\y^D}\right)+\bdy a_{\w^D,\y^D}.
\end{equation}

The first term in Formula~(\ref{eq:3}) is a sum of terms coming from
pairs $(A,B)$ where one of the following cases holds.
\begin{enumerate}
\item $A\in\RectEmpt(\w^D,\x^D)$ and $B\in\RectEmpt(\x^D,\y^D)$ (for some
  $\x^D$). These contributions cancel in pairs exactly as in
  \cite[Proposition 2.8]{MOST07:CombinatorialLink}; see
  Figure~\ref{fig:d-square-zero}.
\item\label{item:RectHalf} $A\in\RectEmpt(\w^D,\x^D)$ and
  $B\in\HalfEmpt(\rho_{i,j};\x^D,\y^D)$ (for some $\x^D$ and
  $\rho_{i,j}$). There are several cases here, the most interesting of
  which is illustrated in Figure~\ref{fig:D-forces}(c).  In this case,
  the relation $\rho_{S,i,j}\cdot\rho_{j,m}=\rho_{S,i,m}$ implies this term
  cancels with a pair of half strips $(A',B')$ obtained by cutting the
  domain horizontally instead of vertically.
\item\label{item:HalfHalf} $A\in\HalfEmpt(\rho_{i,j};\w^D,\x^D)$ and
  $B\in\HalfEmpt(\rho_{l,m};\x^D,\y^D)$ (for some $\x^D$,
  $\rho_{i,j}$, and $\rho_{l,m}$). Again, there are several cases. The
  two half-strips may be disjoint (Figure~\ref{fig:D-forces}(a)), or
  they may form a sideways ``T'' (Figure~\ref{fig:D-forces}(b)); in
  these two cases, relation~\eqref{eq:rel1} implies the contributions
  from taking the two strips in the two different orders cancel. The
  two half-strips may abut top to bottom, in an ``L''-shape
  (Figure~\ref{fig:D-forces}(c)); this cancels with one of the cases
  from Item~(\ref{item:RectHalf}).

  Another possibility is that the upper right corner of $B$ is the
  lower right corner of $A$, as in Figure~\ref{fig:D-forces}(d). This
  configuration contributes a coefficient of
  $\rho_{j,l}\cdot\rho_{i,j}$ (times some $U$-power). There is also a
  half-strip, $A\cup B$, which contributes $\rho_{i,l}$ to $\bdy\w$;
  since $\bdy\rho_{i,l}=\rho_{j,l}\cdot\rho_{i,j}$ in this case, these
  terms cancel.

  Finally, the half strips may overlap as in
  Figure~\ref{fig:D-forces}(e). But in this case the coefficient
  contributed is $\rho_{i,l}\cdot\rho_{j,m}$ which is $0$.
\end{enumerate}

Note that all terms in $\bdy a_{\w^D,\y^D}$ cancelled against terms in
Part~(\ref{item:HalfHalf}). This completes the proof.
\end{proof}

Finally, we turn to the gradings on $\CFPDm(\HD^D)$. Fix any
planar grid diagram $\HD=(\RR^2,\XX,\OO)$ such that
$\HD^D=(H^D,\XX^D,\OO^D)$ can be obtained by cutting $\HD$. Then,
for a generator $\x^D\in\S(\HD^D)$, there are numbers $\cI(\XX,\x^D)$
and $\cI(\OO,\x^D)$, as in Section~\ref{sec:background}. These numbers
obviously do not depend on the choice of $\HD$. Further, fix any
generator $\x\in\S(\HD)$ extending $\x^D$. Then we have a number
$\cI(\x,\x^D)$, which again does not depend on the choice of $\HD$ or
$\x$. Now, define the gradings of $\x^D$ by
\begin{align*}
  A(\x^D)&=\cI(\XX,\x^D)-\cI(\OO,\x^D)\\
  \mu(\x^D)&=\cI(\x,\x^D)-2\cI(\OO,\x^D).
\end{align*}
Extend these definitions to all of $\CFPDm(\HD^D)$ by setting
$A(a\x^D)=A(a)+A(\x^D)$ and $\mu(a\x^D)=\mu(a)+\mu(\x^D)$ for $a\in\Alg_{N,k}$.

\begin{proposition}\label{prop:D-grading}The gradings $A$ and $\mu$ make $\CFPDm(\HD^D)$
  into a graded module over $\Alg_{N,k}$. The differential $\bdy$ on
  $\CFPDm(\HD^D)$ drops $\mu$ by $1$ while preserving $A$.
\end{proposition}
(When assigning gradings to the algebra, we let $L_X$ denote the set
of $i-1/2$ which are \emph{not} $y$-coordinates of points in $\XX^D$,
and similarly for $L_O$.)
\begin{proof}
  The first statement is trivial. To verify that the differential
  drops $\mu$ by $1$, write $\x=(\x^A,\x^D)$. Suppose that
  $\left(\prod_\ell U_\ell^{n_\ell}\right)\rho_{i,j}\cdot\y^D$ occurs in
  $\bdy\x^D$. Then
  \[
  \cI(\x,\x^D)-\cI(\y,\y^D)=1+\#\{(r,s)\in\x^A\mid i< s<j\}.
  \] 
  This is exactly $1+\Cross(\rho_{S,i,j})$, where
  $S=\{0,\dots,N\}\setminus\Image(\sigma_{\x^D})$. Also,
  \[
  \cI(\OO,\x^D)-\cI(\OO,\y^D)=\left(\sum_\ell n_\ell\right)+L_O(\rho_{i,j}).
  \]
  This implies that the differential decreases $\mu$ by $1$, as
  desired.  That the differential preserves $A$ is similar but easier.
\end{proof}

%%% Local Variables: 
%%% mode: latex
%%% TeX-master: "PlanarMain"
%%% End: 

\section{The Type \textalt{$A$}{A} module}
\label{sec:type-A}
The module $\CFPAm$ is much smaller than $\CFPDm$. Fix a partial
planar grid diagram $\HD^A$ with width $k$ and height $N+1$. The
module $\CFPAm(\HD^A)$ is freely generated over $\Ring$ by
$\S(\HD^A)$. There is a differential $\bdy$ on $\CFPAm(\HD^A)$ defined
by
\[
\bdy\x^A=\sum_{\y^A\in\S(\HD^A)}\sum_{
  \begin{subarray}{c}
    R\in\RectEmpt(\x^A,\y^A)\\
    \XX(R)=0
  \end{subarray}
}U(R)\cdot\y^A.
\]
It remains to define an action of $\Alg_{N,k}$ on $\CFPAm(\HD^A)$.

Given a generator $\x^A\in\S(\HD^A)$, let $\sigma_{\x^A}$ denote the
corresponding map $\{0,\dots,k-1\}\allowbreak\to\{0,\dots,N\}$. We define an
action of the idempotents $\Idem_{N,k}$ by
\[
\x^AI_S=
\begin{cases}
  \x^A & \text{if $S=\Image(\sigma_{\x^A})$}\\
  0 & \text{otherwise.}
\end{cases}
\]
This is, in some sense, exactly the opposite of the action
of the idempotents on $\CFPDm$.

\begin{figure}
  \centering
  \includegraphics{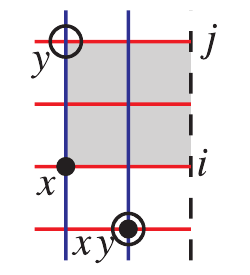}
  \caption{\textbf{An element of $\Half(\x,\y;\rho_{i,j})$.} The
    definition is essentially the same as the definition for $\CFPDm$,
    only rotated by $180$ degrees.}
  \label{fig:AHalfStrip}
\end{figure}

Given generators $\x^A$ and $\y^A$ in $\S(\HD^A)$ and a generator
$\rho_{i,j}$ of $\Alg_{N,k}$ (which we view as a chord in $Z$ from
$y=i$ to $y=j$) define $\Half(\x,\y;\rho_{i,j})$ to be empty unless
$x_k = y_k$ for all but one~$k$, and in this case let it be the singleton
set containing the rectangle (or ``half-strip'') $H$ with lower left
corner $x_k$ and upper left corner $y_k$, and right edge $\rho_{i,j}$
if such a rectangle exists, and empty otherwise. See
Figure~\ref{fig:AHalfStrip}. Call a half strip
$H\in\Half(\x^A,\y^A;\rho_{i,j})$ \emph{empty} if the interior of $H$
is disjoint from $\x^A$. Let $\HalfEmpt(\x^A,\y^A;\rho_{i,j})$ denote
the set of empty half strips in $\Half(\x^A,\y^A;\rho_{i,j})$.

We define an action by the generators $\rho_{i,j}$ of $\Alg_{N,k}$ by
\[
\x^A\rho_{i,j}=\sum_{\begin{subarray}{c}\y^A\in\S(\Sigma)\\
    H\in\HalfEmpt(\x^A,\y^A;\rho_{i,j})\\ \XX(H)=0\end{subarray}}
U(H)\cdot\y^A.
\]
(The sum contains at most one term.)

\begin{proposition}The module $\CFPAm(\HD^A)$ is a differential
  $\Alg_{N,k}$-module. That is:
  \begin{enumerate}
  \item\label{item:resp-rel} The action of the $\rho_{i,j}$ defined
    above respects the relations in $\Alg_{N,k}$.
  \item\label{item:leibniz} The action satisfies the Leibniz rule.
  \item\label{item:bdy-sqr} The differential $\bdy$ satisfies
    $\bdy^2=0$.
  \end{enumerate}
\end{proposition}
\begin{proof}
  (The reader may wish to compare this with the proof of
  Proposition~\ref{prop:D-d-squared}: the pictures are almost the
  same, but their interpretations are different.)

  That the $\Alg_{N,k}$-action respects the three
  relations~\eqref{eq:rel1},~\eqref{eq:rel2} and~\eqref{eq:rel3}
  follow from the cases illustrated in
  Figure~\ref{fig:ARespRel}. In parts (a) and (b), we have
  \[
  \left(\{a,c\}\rho_{i,j}\right)\rho_{l,m}=\left(\{a,c\}\rho_{l,m}\right)\rho_{i,j}=\{b,d\}.
  \] 
  so relation~\eqref{eq:rel1} is respected.
  (We suppress the $U$-powers, but since these depend only on the
  domains they, too, agree.)

  \begin{figure}
    \centering
    %Font is 12 point.
    \includegraphics{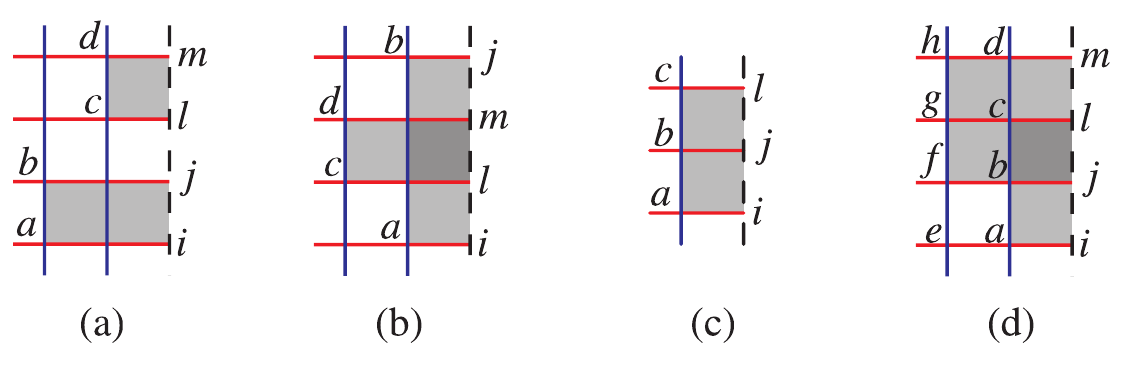}
    \caption{\textbf{The $\Alg_{N,k}$-action on $\CFPAm$ respects the
        relations on the algebra.} Parts (a) and (b) correspond to
      relation~\eqref{eq:rel1}. Part (c) corresponds to
      relation~\eqref{eq:rel2}. Part (d) corresponds to
      relation~\eqref{eq:rel3}.}
    \label{fig:ARespRel}
  \end{figure}
  In part (c) of Figure~\ref{fig:ARespRel},
  \[
  (\{a\}\rho_{i,j})\rho_{j,l}=\{b\}\rho_{j,l}=\{c\}=\{a\}\rho_{i,l},
  \]
  so relation~\eqref{eq:rel2} is respected.

  In part (d) of Figure~\ref{fig:ARespRel} we have
  \[
  \left(\{a,f\}\rho_{i,l}\right)\rho_{j,m}=\{c,f\}\rho_{j,m}=0
  \]
  since the corresponding half-strip is not empty.  So,
  relation~\eqref{eq:rel3} is respected. (This is only one of the two
  pictures we need to check in this case, but the other is similar.)

  This proves Part~(\ref{item:resp-rel}).
  
  Part~(\ref{item:leibniz}) follows from
  Figure~\ref{fig:A-Leibniz}. More precisely, it suffices to show that
  for any $i,j$,
  \[\bdy\left(\x^A\rho_{i,j}\right)=\left(\bdy\x^A\right)\rho_{i,j}+\x^A\left(\bdy\rho_{i,j}\right).\]
  Both $\bdy\left(\x^A\rho_{i,j}\right)$ and
  $\left(\bdy\x^A\right)\rho_{i,j}$ correspond to a domain which is a
  union of a rectangle and a half-strip. The most interesting case is
  when these abut to form an ``L''-shape, as in
  Figure~\ref{fig:A-Leibniz}. There, for $\x^A=\{b,d\}$ we
  have
  \begin{align*}
    \bdy\{b,d\}&=\{a,e\}\\
    \{a,e\}\rho_{i,l}&=\{c,e\}\\
    \{b,d\}\rho_{j,l}\rho_{i,j}&=\{c,e\}\\
    \{b,d\}\rho_{i,l}&=0,
  \end{align*}
  so
  \[
  \bdy\left(\{b,d\}\rho_{i,l}\right)=0=\left(\bdy\{b,d\}\right)\rho_{i,l}
  +\{b,d\}\left(\bdy\rho_{i,l}\right).
  \]
  (The other interesting but similar case is obtained by flipping
  Figure~\ref{fig:A-Leibniz} vertically.) This proves Part~\ref{item:leibniz}.

  \begin{figure}
    \centering
    \includegraphics{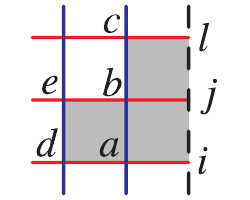}
    \caption{\textbf{The Leibniz rule for $\CFPAm$.} The
      domain shown can be decomposed in two ways: as a rectangle
      followed by a half-strip, or as two half strips. These
      correspond to $\left(\bdy\{b,d\}\right)\rho_{i,l}$ and
      $\{b,d\}\left(\bdy\rho_{i,l}\right)$ respectively. }
    \label{fig:A-Leibniz}
  \end{figure}

  Part~(\ref{item:bdy-sqr}) follows from the same argument as in the
  closed case \cite[Proposition 2.8]{MOST07:CombinatorialLink}; see
  also Figure~\ref{fig:d-square-zero}.
\end{proof}

Finally, we turn to the gradings on $\CFPAm(\HD^A)$. Define
\begin{align*}
A(\x^A)&=\cI(\XX^A,\x^A)-\cI(\OO^A,\x^A)\\
\mu(\x^D)&=\cI(\x^A,\x^A)-2\cI(\OO^A,\x^A).
\end{align*}
\begin{proposition}These gradings make $\CFPAm(\HD^A)$ into a graded
  $\Alg_{N,k}$-module. The differential on $\CFPAm(\HD^A)$ preserves
  the Alexander grading $A$ and drops the Maslov grading $\mu$ by $1$.
\end{proposition}
(When assigning gradings to the algebra, we let $L_X$ denote the set
of $i-1/2$ which are $y$-coordinates of points in $\XX^A$, and
similarly for $L_O$.)
\begin{proof}
  We check that multiplication preserves the $A$ grading. Suppose
  $\x^A\rho_{i,j}=\left(\prod_\ell U_\ell^{n_\ell}\right)\y^A$. Then
  \begin{align*}
    \cI(\XX^A,\x^A)&=\cI(\XX^A,\y^A)-L_X(\rho_{i,j})\\
    \cI(\OO^A,\x^A)&=\cI(\OO^A,\y^A)-L_O(\rho_{i,j})+\sum_\ell n_\ell.
  \end{align*}
  The result follows.

  That multiplication preserves $\mu$ is similar; see also the proof
  of Proposition~\ref{prop:D-grading} That the differential preserves
  $A$ and drops $\mu$ by $1$ is straightforward.
\end{proof}

\begin{remark}The definition of $\CFPAm$ is somewhat different in
  spirit from the definition of $\CFAa$ for bordered three-manifolds
  in \cite[Section 7]{LOT1}: there the product $\x^A a$ is defined
  directly for any algebra element $a$. In our setting, we could do
  this by counting more complicated domains than rectangles.
\end{remark}

%%% Local Variables: 
%%% mode: latex
%%% TeX-master: "PlanarMain"
%%% End: 

\section{The pairing theorem}
\begin{theorem}Let $\HD$ be a planar grid diagram, decomposed as
  $\HD^A\cup_Z\HD^D$, where $\HD^A$ (respectively $\HD^D$) is a
  partial planar grid diagram with width $k$ (respectively $N+1-k$)
  and height $N+1$. Then
  \[
  \CFPm(\HD)\cong\CFPAm(\HD^A)\otimes_{\Alg_{N,k}}\CFPDm(\HD^D),
  \]
  as $(\ZZ\oplus \ZZ)$-graded chain complexes over $\Ring$.
\end{theorem}
\begin{proof}
  There is an obvious identification between the generators of
  $\CFPm(\HD)$ and the generators of
  $\CFPAm(\HD^A)\otimes_{\Alg_{N,k}}\CFPDm(\HD^D)$. It follows from
  their definitions that this identification respects the $A$ and
  $\mu$ gradings.

  The rest of the proof is essentially trivial, so we write it with
  formulas to make it seem complicated. Given a generator $\x$ of
  $\CFPm(\HD)$, we split $\bdy\x$ into three pieces, according to
  whether the domain rectangle is entirely to the left of the dividing
  line $\{x = k-1/4\}$, crosses the dividing line, or is entirely to
  the right of the dividing line:
  \[
  \bdy\x = \bdy_L\x + \bdy_M\x + \bdy_R\x.
  \]
  Then if $\x$ is identified with $\x^A\otimes\x^D$, we have
  \begin{align*}
    \bdy(\x^A\otimes\x^D)&=(\bdy\x^A)\otimes\x^D+\x^A\otimes(\bdy\x^D)\\
    &= \bdy_L\x + \bdy_R\x
    +\sum_{\y^D\in\S(\HD^D)}\sum_{
      \begin{subarray}{c}
        H\in\HalfEmpt(\rho_{i,j};\x^D,\y^D)\\
        \XX(H)=0
      \end{subarray}
}U(H)(\x^A\otimes
    \rho_{i,j}\y^D)\\
    &= \bdy_L\x + \bdy_R\x
    +\sum_{\y^D\in\S(\HD^D)}\sum_{
      \begin{subarray}{c}
        H\in\HalfEmpt(\rho_{i,j};\x^D,\y^D)\\
        \XX(H)=0
      \end{subarray}
}U(H)(\x^A\rho_{i,j}\otimes\y^D)\\
    &=\bdy_L\x+\bdy_R\x\\
    &\qquad+\sum_{
       \substack{
         \y^D\in\S(\HD^D)\\
         \y^A\in\S(\HD^A)\\
         \rho_{i,j}
       }
    }\,\sum_{
      \begin{subarray}{c}
        H\in\HalfEmpt(\rho_{i,j};\x^D,\y^D)\\
        \XX(H)=0
      \end{subarray}
}\,\sum_{
  \begin{subarray}{c}
    H'\in\HalfEmpt(\x^A,\y^A;\rho_{i,j})\\
    \XX(H')=0
  \end{subarray}
}U(H'\cup H)(\y^A\otimes
      \y^D)\\
      &=\bdy_L\x + \bdy_R\x + \bdy_M\x\\
      &=\bdy\x,
  \end{align*}
  as desired.
\end{proof}

\begin{remark}More useful is the fact that $\CFPm(\HD)$ is
  quasi-isomorphic to the derived tensor product
  $\CFPAm(\HD^A)\DTP_{\Alg_{N,k}}\CFPDm(\HD^D)$. For instance, this
  allows one to simplify the complexes $\CFPAm$ and $\CFPDm$ more
  dramatically before taking the tensor product.
In fact, the $\Alg_{N,k}$-module $\CFPDm(\HD^D)$ is projective (hence
flat), as one can see by imitating an argument from Bernstein and
Lunts \cite[Proposition 10.12.2.6]{BernsteinLunts94:EquivariantSheaves}.
It follows that the derived tensor product agrees with the ordinary one.
\end{remark}
%%% Local Variables: 
%%% mode: latex
%%% TeX-master: "PlanarMain"
%%% End: 

\section{Bimodules}
At this point we have encountered left and right modules over
$\Alg_{N,k}$. We will now see that bimodules also have several
important roles to play. (The material in this section is analogous to
material in~\cite{LOT2}.)

\subsection{Freezing}
\label{sec:freezing}
We have studied how to take a planar grid diagram and make a single
vertical cut. In the spirit of factoring a braid into generators,
however, we might want to make several different vertical cuts. In
this section we will see that the correct objects to assign to slices
in the middle are $(\Alg_{N,k}$,$\Alg_{N,l})$-bimodules.

That is, consider the result of slicing a planar grid diagram $\HD$ along
the lines $Z_1=\{x=k-1/4\}$ and $Z_2=\{x=l-1/4\}$ (with $l>k$). The
result is two partial planar grid diagrams $\HD^A=\HD\cap\{x<k-1/4\}$
and $\HD^D=\HD\cap\{x>l-1/4\}$, and a \emph{middle partial planar grid
  diagram} $\HD^{DA}=\HD\cap\{k-1/4<x<l-1/4\}$. We will associate an
$(\Alg_{N,k},\Alg_{N,l})$-bimodule $\CFPDAm(\HD^{DA})$ to $\HD^{DA}$.

A generator for $\HD^{DA}$ is an $(l-k)$-tuple of points
$\x=\{x_i\}_{i=k}^{l-1}$; a generator $\x$ corresponds to an
injection $\sigma_{\x}\co\{k,\dots,l-1\}\to\{1,\dots,N\}$. (For
consistency with earlier notation, we should really write $\x$ as
$\x^{DA}$, but the notation becomes too cumbersome.) Let
$\S(\HD^{DA})$ denote the set of generators for $\HD^{DA}$. Call a
generator $\x$ \emph{compatible} with an idempotent
$I_S\in\Alg_{N,k}$ if $\Image(\sigma_{\x})\cap S=\emptyset$. As a
left module, $\CFPDAm(\HD^{DA})$ is a direct sum of elementary
modules,
\[
\CFPDAm(\HD^{DA})=\bigoplus_{\begin{subarray}{c}\x\in\S(\HD^{DA})\\
    S\text{ compatible with }\x \end{subarray}}\Alg_{N,k}I_S.
\]
We will write the generator of the summand $\Alg_{N,k}I_S$ coming from
$\x$ as $I_S\x.$ Note that, unlike for $\CFPDm$ or $\CFPAm$,
the generator $\x$ does not determine the idempotent $S$.

We next define a differential on $\CFPDAm(\HD^{DA})$. Given generators
$\x,\y\in\S(\HD^{DA})$ such that $x_n=y_n$ for $n\neq m$ (for some $m$),
and $i<j\in\{0,\dots,N\}$, define $\Half(\rho_{i,j};\x,\y)$ to be the
set of rectangles with upper right corner at $x_m$, lower right corner
at $y_m$ and left edge the segment $\rho_{i,j}$ in $Z_1$ from
$(k-1/4,i)$ to $(k-1/4,j)$. Define $\HalfEmpt(\rho_{i,j};\x,\y)$ to be
the subset of $\Half(\rho_{i,j};\x,\y)$ consisting of empty
half-strips, i.e., half strips not containing any element of $\x$ in
their interiors. Then set
\[
\bdy (I_S\x)=
\sum_{\y\in\S(\HD^{DA})}\,
\sum_{
  \begin{subarray}{c}
    R\in\RectEmpt(\x,\y)\\
    \XX(R)=0
  \end{subarray}}
  U(R)\cdot I_S\y
+
\sum_{
  \begin{subarray}{c}
    \y\in\S(\HD^{DA})\\
    i<j\in\{0,\dots,N\}
  \end{subarray}
}
\sum_{
  \begin{subarray}{c}
    H\in\HalfEmpt(\rho_{i,j};\x,\y)\\
    \XX(H)=0
  \end{subarray}}
  U(H)\cdot I_S\rho_{i,j}\y.
\]
Here, the notation $I_S\rho_{i,j}\y$, though suggestive, should
be explained. If $i\in S$ and $j\notin S$ then $I_S\rho_{i,j}\y$
denotes $\rho_{i,j}I_{T}$, where $T=(S\setminus i)\cup j$, if $T$ is
compatible with $\y$. Otherwise (i.e., if $i\notin S$, $j\in S$,
or $T$ is not compatible with $\y$) we declare
$I_S\rho_{i,j}\y$ to be $0$.

Finally, we define the right module structure on
$\CFPDAm(\HD^{DA})$. Given a primitive idempotent $I_T\in\Idem_{N,l}$,
define
\[
(I_S\x)I_T =
\begin{cases}
  I_S\x & \text{if $\left(S\cup
      \Image(\sigma_{\x})\right)\cap T=\emptyset$}\\
  0 & \text{otherwise}.
\end{cases}
\]
Given generators
$\x,\y\in\S(\HD^{DA})$ such that $x_\ell=y_\ell$ for $\ell\neq m$
(for some $m$),
and $i<j\in\{0,\dots,N\}$, define $\Half(\x,\y;\rho_{i,j})$ to be the
set of rectangles with lower left corner at $x_m$, upper left corner
at $y_m$ and right edge the segment $\rho_{i,j}$ in $Z_2$ from
$(l-1/4,i)$ to $(l-1/4,j)$. Define $\HalfEmpt(\x,\y;\rho_{i,j})$ to be
the subset of $\Half(\x,\y;\rho_{i,j})$ consisting of empty
half-strips, i.e., half strips not containing any element of $\x$ in
their interiors.

Given a chord $\rho_{i,j}$ in $Z_2$ from $(l-1/4,i)$ to $(l-1/4,j)$,
define $\Strip(\rho_{i,j})$ to be the horizontal strip with right
edge $\rho_{i,j}\subset Z_2$ and left edge $\rho_{i,j}\subset
Z_1$. Given $\rho_{i,j}$ and a generator $\x\in\S(\HD^{DA})$,
define $\StripEmpt(\x;\rho_{i,j})$ to be the empty set if
$\Strip(\rho_{i,j})$ contains a point in $\x$ (even along its
boundary) and the singleton set $\Strip(\rho_{i,j})$ if
$\Strip(\rho_{i,j})$ does not contain a point in $\x$.

At last, define
\[
(I_S\x)\rho_{i,j}=
\sum_{
  \begin{subarray}{c}
    E\in\StripEmpt(\x,\rho_{i,j})\\
    \XX(E)=0
  \end{subarray}
}
U(E)\cdot I_S\rho_{i,j}\x
+
\sum_{\y\in\S(\HD^{DA})}
\sum_{
  \begin{subarray}{c}
    H\in\HalfEmpt(\x,\y;\rho_{i,j})\\
    \XX(H)=0
  \end{subarray}}
  U(H)\cdot I_S\y.
\]

These definitions are, in fact, compatible:
\begin{proposition}The module $\CFPDAm(\HD^{DA})$ is a differential
  $(\Alg_{N,k},\Alg_{N,l})$-bimodule.
\end{proposition}
We leave the proof to the interested reader.

As on $\CFPDm$, the grading of a generator $I_S\x$ of
$\CFPDAm(\HD^{DA})$ is given by
\begin{align*}
  A(I_S\x)&=\cI(\XX,\x)-\cI(\OO,\x)\\
  \mu(I_S\x)&=\cI(\bar{\x},\x)-2\cI(\OO,\x),
\end{align*}
where $\HD=(\RR^2,\XX,\OO)$ is any planar diagram completing $\HD^{DA}$,
and $\bar{\x}$ is a generator in $\S(\HD)$ completing $\x$ and
compatible with the idempotent $I_S$ in the obvious sense.

Finally, the module $\CFPDAm(\HD^{DA})$ satisfies a pairing
theorem:
\begin{proposition} With notation as above,
  \begin{align*}
    \CFPAm(\HD^A\cup_{Z_1}\HD^{DA})&=\CFPAm(\HD^A)\otimes_{\Alg_{N,k}}\CFPDAm(\HD^{DA})\\
    \CFPDm(\HD^{DA}\cup_{Z_2}\HD^D)&=\CFPDAm(\HD^{DA})\otimes_{\Alg_{N,l}}\CFPDm(\HD^D)\\
    \CFPm(\HD)&=\CFPAm(\HD^A)\otimes_{\Alg_{N,k}}\CFPDAm(\HD^{DA})\otimes_{\Alg_{N,l}}\CFPDm(\HD^D).
  \end{align*}
\end{proposition}
The proof is obvious. The analogous result for cutting along more than
two vertical lines is also true.

\begin{remark}The notation $\CFPDAm$ denotes that the module is ``Type
  $D$'' from the left and ``Type $A$'' from the right.
\end{remark}

\subsection{Type \textalt{$A$}{A} to Type \textalt{$D$}{D}}
\label{sec:A-to-D}
The reader might wonder about the relation between $\CFPAm$ and
$\CFPDm$. One might expect that they are, in some appropriate
sense, dual to each other. In the case of bordered Heegaard Floer
homology this is true. In this section, we hint at that story by
reconstructing $\CFPDm$ from $\CFPAm$. In Section~\ref{sec:D-to-A} we
will discuss going the other direction. We will suppress both the
gradings and the $U$-variables: our treatment of both has been too
na\"ive to extend properly to the present discussion.

Let $k'=N+1-k$.
We construct a $(\Alg_{N,k},\Alg_{N,k'})$-bimodule $\CFPDDm_{N,k}$
so that
$\CFPDm(\HD^D)=\CFPAm(\HD^D)\otimes_{\Alg_{N,k}}\CFPDDm_{N,k}$. Actually,
unlike a traditional bimodule with a left action and a right action,
we will construct the $\Alg_{N,k}$- and $\Alg_{N,k'}$-actions as a
pair of commuting \emph{left} actions, so the module
$\CFPAm(\HD^D)\otimes_{\Alg_{N,k}}\CFPDDm_{N,k}$ comes equipped with a
left action rather than a right action.

\begin{figure}
  \centering
  %Font is 12 point.
  \includegraphics{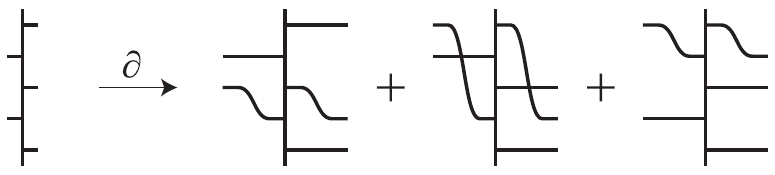}
  \caption{\textbf{A graphical representation of $\CFPDDm_{N,k}$.} The
    case shown is $N=4$, $k=2$. The element $I_S$, for $S=\{1,3\}$, is
    shown on the left. On the right is the differential of $I_S\otimes
    1$. This
    graphical representation treats $\CFPDDm_{N,k}$ as a traditional
    (left,right) bimodule, rather than a (left,left) bimodule; this is
    the reason that the strands on the right are downward-veering.}
  \label{fig:CPDD}
\end{figure}

The module $\CFPDDm_{N,k}$ is easy to describe. Note that there is an
obvious isomorphism $\Idem_{N,k}\to\Idem_{N,k'}$, taking $I_S$ to
$I_{\{0,\dots,N\}\setminus S}$. This makes $\Alg_{N,k'}$ into a
right $\Idem_{N,k}$-module. The module $\CFPDDm_{N,k}$ is just
\[
\Alg_{N,k'}\otimes_{\Idem_{N,k}}\Alg_{N,k'}
\]
where the tensor product identifies the \emph{right} actions of
$\Idem_{N,k}$ on $\Alg_{N,k'}$ and $\Alg_{N,k'}$. This module,
then, is equipped with two left actions. The differential on
$\CFPDDm_{N,k}$ is not the one inherited from the tensor
product. Rather, for $S$ a $k$-element subset of $\{0,\dots,N\}$ we
define
\[
\bdy (I_S\otimes 1)=\sum_{\begin{subarray}{c}i\in S\\ j>i\\ j\notin
    S\end{subarray}}\rho_{i,j}^{k}\rho_{i,j}^{k'} (I_T\otimes 1)
\]
where $\rho_{i,j}^k$ denotes the element $\rho_{i,j}$ of $\Alg_{N,k}$,
$\rho_{i,j}^{k'}$ denotes the element $\rho_{i,j}$ of
$\Alg_{N,k'}$, and $T=(S\setminus i)\cup j$. We extend the
differential to all of $\CFPDDm_{N,k}$ by the Leibniz rule. An example
is illustrated in Figure~\ref{fig:CPDD}.

\begin{lemma}The module $\CFPDDm_{N,k}$ is a differential
  $(\Alg_{N,k},\Alg_{N,k'})$-bimodule.
\end{lemma}
\begin{proof}
This is immediate from the definitions.
\end{proof}

One can view the module $\CFPDDm_{N,k}$ as the (Type $D$, Type
$D$) module associated to a middle partial planar grid diagram with
zero $\beta$-lines (i.e., in the notation of
Section~\ref{sec:freezing}, $k=l$). The generator corresponds to the
empty set in $\alphas\cap\betas$. The differential comes from the
strips $\Strip(\rho_{i,j})$ (as in Section~\ref{sec:freezing}).

As promised, we have the following pairing theorem:
\begin{proposition}Fix a partial Heegaard diagram $\HD^D$. Then
  \[\CFPDm(\HD^D)=\CFPAm(\HD^D)\otimes_{\Alg_{N,k}}\CFPDDm_{N,k}.\]
\end{proposition}
\begin{proof}
  The tensor product $\CFPAm(\HD^D)\otimes_{\Alg_{N,k}}\CFPDDm_{N,k}$
  is a direct sum of elementary modules $\x^A\otimes\Alg_{N,k'}I_S$, one for
  each generator $\x^A$ of $\CFPAm(\HD^D)$, where
  $S=\{0,\dots,N\}\setminus\Image(\sigma_{\x^A})$. The part of the
  differential on the tensor product coming from the differential on
  $\CFPAm(\HD^D)$ counts empty rectangles. The part of the
  differential coming from the differential on $\CFPDDm_{N,k}$ counts
  empty half strips, exactly as on $\CFPDm(\HD^D)$.
\end{proof}

\subsection{Remarks on Type \textalt{$D$}{D} to Type \textalt{$A$}{A}}
\label{sec:D-to-A}
Turning the module $\CFPDm$ into $\CFPAm$ is more subtle than turning
$\CFPAm$ into $\CFPDm$. It is clearly not possible to find a module
$\CFPAAm_{N,k}$ so that $\CFPAm$ is exactly equal to
$\CFPAAm_{N,k}\otimes_{\Alg_{N,k'}}\CFPDm$: the ranks of these modules
over $\Ring$ prevent this.

There are two approaches one might take. One approach is to use
properties of $\CFPDDm_{N,k}$ to prove that it induces an equivalence
of categories
$\DerBounded(\Alg_{N,k}-\ModCat)\to\DerBounded(\Alg_{N,k'}-\ModCat)$,
and then construct a bimodule giving the inverse equivalence of
categories.
  
Another approach is to define $\CFPAAm_{N,k}$ as an
$\Ainf$-bimodule. Using the appropriate model for the $\Ainf$-tensor
product (see, e.g., \cite[Section 2]{LOT1}), it is then possible for
$\CFPAm$ to be exactly
$\CFPAAm_{N,k}\otimes_{\Alg_{N,k'}}\CFPDm$. The generators and first few
$\Ainf$-operations for this $\CFPAAm_{N,k}$ are easy to guess. As an
$\Ring$-module, $\CFPAAm_{N,k}$ would be just
$\Idem_{N,k}\cong\Idem_{N,k'}$. The first few $\Ainf$-relations would be
\begin{align*}
m_1(I_S) &= 0\\
m_2(a,I_S) &=0\\
m_2(I_S, a) &=0\\
m_3(\rho_{i,j},I_S,\rho_{i,j})&=
\begin{cases}
  I_T & \text{ if $i\in S$, $j\notin S$, and where $T=(S\setminus i) \cup j$}\\
  0 & \text{otherwise.}
\end{cases}
\end{align*}
(Even though $\CFPAAm_{N,k}$ should really have two right actions, for
clarity we have written it with one right and one left action.)

Unfortunately, higher $\Ainf$-relations are harder to guess and, at
least in the case of bordered Heegaard Floer homology, depend on some
choices. Fortunately, in the case of bordered Heegaard Floer homology,
these modules are induced by counts of holomorphic curves, so we need
not build them by hand; see~\cite{LOT1}. (In particular, it turns out
that the choices are induced by a choice of almost complex structure.)
The challenge in defining $\CFPAAm_{N,k}$, then, becomes counting
holomorphic curves.

%%% Local Variables: 
%%% mode: latex
%%% TeX-master: "PlanarMain"
%%% End: 

\section{How the real world is harder}
In this section, we preview the difficulties involved in using the
ideas from this paper to define more useful invariants.

\subsection{Complications for \textalt{$\HFa$}{HF-hat} of
  \textalt{$3$}{3}-manifolds}
As discussed in the introduction, applying the ideas of this paper to
the case of the Heegaard Floer group $\HFa(Y)$ gives an invariant of
$3$-manifolds with boundary; see~\cite{LOT1}. The main complications
are as follows.
\subsubsection{Heegaard diagrams.} Instead of working with grid
diagrams, the invariant $\HFa(Y^3)$ is defined by using a ``Heegaard
diagram'' for $Y$. One needs, then, an appropriate family of partial
Heegaard diagrams. Such a class, called either ``Heegaard diagrams
with boundary'' or ``bordered Heegaard diagrams'' was presented
in~\cite{Lipshitz06:BorderedHF};  see also~\cite[Section
4]{LOT1}. These diagrams are induced by a self-indexing Morse function
$f$ on a three manifold with boundary $(Y,\bdy Y)$ such that $\nabla
f$ is tangent to $\bdy Y$ (and subject to a few more
constraints). Bordered Heegaard diagrams specify not just the
three-manifold $Y$ but also a parametrization of $\bdy Y$; this is
obviously needed for the pairing theorem to make sense.

One incidental effect is that the algebra $\Alg_{N,k}$ needs to be
modified somewhat. In the planar setting, each $\alpha$-line
intersects the interface $Z$ in a single point; in the bordered case
(or the toroidal case) this is not true. The solution in the bordered
case is to work with a subalgebra of $\Alg_{N,k}$ which, roughly,
remembers how the points $\alphas\cap Z$ are paired-up. (In the
toroidal case described below, it is more convenient to remember only
half of the
points and drop the requirement that strand diagrams be
upward-veering.)

\subsubsection{Holomorphic curves.} Like the closed Heegaard Floer
invariant $\CFa(Y)$, the definitions of the bordered Heegaard Floer
invariants $\CFAa(Y)$ and $\CFDa(Y)$ involve counting holomorphic
curves. The analytic setup here is somewhat nonstandard, complicating
matters.

Like $\CFa(Y)$, the techniques of Sarkar and
Wang~\cite{SarkarWang07:ComputingHFhat} allow one to compute
$\CFAa(Y)$ and $\CFDa(Y)$ combinatorially, by using a particular kind
of diagram called a \emph{nice diagram.} Such diagrams also make the
pairing theorem as trivial as it was in the planar case. However,
there is currently no way to prove invariance for even the closed
invariant while staying in the class of nice diagrams; also, working
with a nice diagram seems to require super-exponentially more
generators in most cases.

\subsubsection{\textalt{$\Ainf$}{A-infinity}-structures and
  noncommutative gradings}
For general Heegaard diagrams, associativity fails for
$\CFAa(Y)$. Fortunately, associativity holds up to homotopy, and in
fact one can organize the higher associators neatly into the structure
of an $\Ainf$-module. (In the case that the bordered Heegaard diagram is
nice, all higher associators vanish, and hence $\CFAa(Y)$ is
an honest module.)

Another algebraic complication is the grading. For boundary of genus
at least one, the algebra $\Alg(F)$ associated to a surface $F$ is
not $\ZZ$-graded but rather is graded by a certain noncommutative
group $G$. (This grading intertwines the homological and $\SpinC$
gradings.) The modules associated to bordered $3$-manifolds are graded by
$G$-sets.

\subsection{Complications for toroidal grid diagrams}
One can also try to pursue an analogue of this theory for toroidal
grid diagrams. Slicing a toroidal grid diagram yields a
representation of a tangle, so this can be viewed as a theory of tangles.
There seem to be two main complications, the second
more serious than the first.
\subsubsection{Boundary degenerations and matrix factorizations.}
For planar grid diagrams, or for bordered Heegaard diagrams, there are
no domains with boundary contained entirely in the $\alpha$-curves (or
entirely in the $\beta$-curves). This prevents certain degenerations
of holomorphic curves (called ``boundary degenerations''
in~\cite{OS05:HFL}). For toroidal grid diagrams, there are such
degenerations. Their cancellation, holomorphically~\cite{OS05:HFL}
or combinatorially~\cite{MOST07:CombinatorialLink}, is delicate,
and not preserved by the slicing operation. The result is that the
invariants one must associate to partial toroidal grid diagrams are
not differential modules but instead matrix factorizations. (Matrix
factorizations also arise in other knot homology theories; see,
e.g.,~\cite{KhovanovRozansky08:MatrixFactorizations}.) Equivalently,
one can deform a suitable version of the algebra $\Alg_{N,k}$ to an
$\Ainf$-algebra with a
nontrivial $\mu_0$.

\subsubsection{Derived equivalences.}
In this paper, we have not talked at all about invariance, because the
planar Floer homology $\CFPm$ is itself not an invariant. For the
toroidal theory, a partial diagram of height $N$ and width $k$ will
result in a module over an algebra $\Alg_{N,k}^{\XX,\OO}$, a variant of
$\Alg_{N,k}$. One can have diagrams for a tangle with different
heights and widths; the ``invariants'' associated to them, then, are
modules over different algebras. In order to even express invariance,
then, one would like derived equivalences
\[
\DerBounded(\Alg_{N,k}^{\XX,\OO}-\ModCat)\to
\DerBounded(\Alg_{N',k'}^{\XX',\OO'}-\ModCat)
\]
between certain of these algebras. Moreover, these must be compatible
with how stabilization acts on the modules. We return to these issues
in a future paper~\cite{LOT2}.

%%% Local Variables: 
%%% mode: latex
%%% TeX-master: "PlanarMain"
%%% End: 

\bibliographystyle{hamsplain}\bibliography{heegaardfloer}
\end{document}